\documentclass[11pt,fleqn,twoside]{article}
\usepackage{amsfonts,amssymb,latexsym}
\makeatletter
\newcommand{\prava}[1]{\small\it
\begin{flushleft}
Copyright \copyright \ 1999 by  #1
\end{flushleft}}

\newcommand{\name}[1]{\begin{flushleft}
                       \LARGE \bf #1
                       \end{flushleft}\vspace{-3mm}}

\newcommand{\Author}[1]{\begin{flushleft}
                       \it #1 \end{flushleft}}

\newcommand{\Adress}[1]{\begin{flushleft}
                       \it #1 \end{flushleft}}

\newcommand{\Date}[1]{\begin{flushleft}
                      \small  \it #1 \end{flushleft}}

\newcommand{\ehkol}{Author \ name}
\newcommand{\ohkol}{Article \ name}
\renewcommand{\@evenhead}{
\hspace*{-3pt}\raisebox{-15pt}[\headheight][0pt]{\vbox{\hbox to \textwidth 
{\thepage \hfil \ehkol}\vskip4pt \hrule}}}
\renewcommand{\@oddhead}{
\hspace*{-3pt}\raisebox{-15pt}[\headheight][0pt]{\vbox{\hbox to \textwidth 
{\ohkol \hfil \thepage}\vskip4pt\hrule}}}
\renewcommand{\@evenfoot}{}
\renewcommand{\@oddfoot}{}

\newcommand{\be}{\begin{equation}}
\newcommand{\ee}{\end{equation}}
\newcommand{\ba}{\hspace*{-5pt}\begin{array}}
\newcommand{\ea}{\end{array}}

\newcommand{\ds}{\displaystyle}
\makeatother

\begin{document}
\thispagestyle{empty}
\setcounter{page}{384}
\renewcommand{\ehkol}{A.K. Prykarpatsky}
\renewcommand{\ohkol}{The Nonabelian Liouville-Arnold Integrability by
Quadratures Problem}

\begin{flushleft}
\footnotesize \sf Journal of Nonlinear Mathematical Physics \qquad
1999, V.6, N~4, \pageref{prykarp-fp}--\pageref{prykarp-lp}.
\hfill {\sc Article}
\end{flushleft}

\vspace{-5mm}

\renewcommand{\footnoterule}{}
{\renewcommand{\thefootnote}{}
 \footnote{\prava{A.K. Prykarpatsky}}}

\name{The Nonabelian Liouville-Arnold Integrability by
Quadratures Problem: a Symplectic Approach}\label{prykarp-fp}

\Author{Anatoliy K. PRYKARPATSKY}

\Adress{Department of Applied Mathematics at the AGH,
Krak\'{o}w 30-059, Poland\\
and Department of Physics at the Eastern Mediterranean University,\\
Gamagusa, Northen Cyprus
\\[1mm]
E-mails: tolik@mailcity.com, prykanat@cybergal.com,
prika@mat.agh.edu.pl, ptykat.as@mozart.emu.edu.tr}

\Date{Received May 10, 1999; Revised: June 23, 1999; Accepted July 22, 1999}

\begin{abstract}
\noindent
A symplectic theory approach is
devised for solving the problem of algebraic-analytical con\-struction of
integral submanifold imbeddings for integrable
(via the nonabelian Liouville-Arnold theorem) Hamiltonian systems
on canonically symplectic phase spa\-ces.
\end{abstract}

\renewcommand{\theequation}{\thesection.\arabic{equation}}
\setcounter{equation}{0}
\setcounter{section}{0}

\section*{0. Introduction}

{\bf 0.1.} As is well known [1, 4], the integrability by quadratures of a
dif\/ferential equation in space ${\mathbb R}^{n}$ is a method of seeking
its solutions by means of f\/inite number of algebraic operations (together
with inversion of functions) and ``quadratures'' -- calculations of integrals
of known functions.

Assume that our dif\/ferential equation is given as a Hamiltonian dynamical
system on some appropriate symplectic manifold  $\left(M^{2n},
\omega ^{(2)}\right)$,  $n\in {\mathbb Z}_{+},$ in the form
\begin{equation}
du/dt=\{H,u\},
\end{equation}
where $u\in M^{2n}$,  $H:M^{2n}\rightarrow {\mathbb R}$ is a
suf\/f\/iciently smooth Hamiltonian function [1, 4] with respect to the Poisson
bracket $\{\cdot ,\cdot \}$  on ${\mathcal D}(M^{2n}),$ dual to the symplectic
structure $\omega ^{(2)}\in \Lambda ^{2}(M^{2n}),$ and
$t\in {\mathbb R}$ is the evolution parameter.

More than one hundred and f\/ifty years ago French mathematicians
and physicists, f\/irst E.~Bour and next J.~Liouville, proved the
f\/irst ``integrability by quadratures'' theorem which in modern
terms~[33] can be formulated as follows.

\medskip

\noindent
{\bf Theorem 0.1.} {\it Let $M^{2n}\simeq T^{\ast }({\mathbb R}^{n})$
be a canonically symplectic phase space and there be given a
dynamical system~(0.1) with a Hamiltonian function $H$: $M^{2n}\times
{\mathbb R}_{t}\rightarrow {\mathbb R}$, possessing a
Poissonian Lie algebra ${\mathcal G}$ of $n \in {\mathbb Z}_{+}$
invariants $H_{j}:M^{2n}\times {\mathbb R}_{t}\rightarrow
{\mathbb R}$, $j=\overline{1,n}$, such that
\begin{equation}
\{H_{i},H_{j}\}=\sum_{s=1}^{n}c_{ij}^{s}H_{s},
\end{equation}
and for all $i,j,k =\overline{1,n}$ the $c_{ij}^{s}\in {\mathbb R}$
are constants on $M^{2n}\times {\mathbb R}_{t}$. Suppose further that
\begin{equation}
M_{h}^{n+1}:=\left\{(u,t)\in M\times {\mathbb R}_{t}:h(H_{j})=h_{j},
\  j=\overline{1,n}, \ h\in {\mathcal G}^{\ast }\right\},
\end{equation}
the integral submanifold of the set ${\mathcal G}$  of invariants
at  a regular element $h \in {\mathcal G}^{\ast }$, is a well defined
connected submanifold of $M\times {\mathbb R}_{t}$. Then, if:

i) all functions of ${\mathcal G}$ are functionally independent on
$M_{h}^{n+1}$;

ii) $\sum\limits_{s=1}^{n}c_{ij}^{s}h_{s}=0$ for all $i,j =\overline{1,n}$;

iii) the Lie algebra ${\mathcal G} = \mbox{\rm span}_{\mathbb R}
\left\{H_{j}:M^{2n}\times {\mathbb R}_{t}\rightarrow {\mathbb R}:
j=\overline{1,n}\,\right\}$ is solvable,
the Hamiltonian system (0.1) on $M^{2n}$  is integrable by quadratures.}

\medskip

As a simple corollary of the Bour-Liouville theorem  one gets the following:

\medskip

\noindent
{\bf Corollary 0.2.} {\it If a Hamiltonian system on $M^{2n}=T^{\ast }
({\mathbb R}^{n})$ possesses just $n\in {\mathbb Z}_{+}$
func\-tionally independent invariants in involution, that is a
Lie algebra ${\mathcal G}$ is abelian, then it is integrable by
quadratures.}

\medskip

In the autonomous case when a Hamiltonian $H=H_{1}$,  and invariants $H_{j}
: M^{2n}\rightarrow {\mathbb R}$,  $j=\overline{1,n}$, are independent of the
evolution parameter $t\in {\mathbb R}$,  the involutivity
condition $\{H_{i},H_{j}\} =0$,  $i,j=\overline{1,n}$,  can be replaced
by the weaker one $\{H,H_{j}\}=c_{j}H$  for some constants $c_{j}\in
{\mathbb R}$, $j=\overline{1,n}$.

The f\/irst proof of Theorem 0.1. was based on a result of S.~Lie, which can
be formulated as follows.

\medskip

\noindent
{\bf Theorem 0.3 (S.~Lie).} {\it  Let vector fields $K_{j}\in \Gamma
(M^{2n})$, $j=\overline{1,n}$, be independent in some open
neighborhood $U_{h}\in M^{2n}$, generate  a solvable Lie algebra
${\mathcal G}$ with respect  to the usual commutator
$[\cdot ,\cdot ]$  on $\Gamma (M^{2n})$ and $[K_{j},K] = c_{j}K$
for all $j=\overline{1,n}$, where $c_{j}\in {\mathbb R}$,
$j=\overline{1,n}$, are  constants. Then the dynamical system
\renewcommand{\theequation}{\thesection.\arabic{equation}${}'$}
\setcounter{equation}{0}
\begin{equation}
du/dt=K(u),
\end{equation}
where $u \in U_{h}\subset M^{2n}$, is
integrable by quadratures. }

\medskip

\noindent
{\bf Example 0.4.} {\bfseries \itshape Motion of three particles
on line ${\mathbb R}$ under uniform potential f\/ield.}
The motion of three particles on the axis ${\mathbb R}$ pairwise  interacting
via a uniform potential f\/ield $Q(\left\| \cdot \right\| )$ is described as
a Hamiltonian system on the canonically symplectic phase space $M =
T^{\ast }({\mathbb R}^{3}$) with the following Lie algebra ${\mathcal G}$ of
invariants on $M^{2n}$:
\renewcommand{\theequation}{\thesection.\arabic{equation}}
\setcounter{equation}{3}
\begin{equation}
\ds H=H_{1}=\sum_{j=1}^{3}p_{j}^{2}/2m_{j}+\sum_{i<j=1}^{3}Q(\left\|
q_{i}-q_{j}\right\| ),
\quad  H_{2}=\sum_{j=1}^{3}q_{j}p_{j},\quad
H_{3}=\sum_{j=1}^{3}p_{j},
\ee
where $(q_{j},p_{j})\in T^{\ast }({\mathbb R})$, $j=\overline{1,3},$ are
coordinates and momenta of particles on the axis ${\mathbb R}$. The
commutation relations for the Lie algebra ${\cal G}$ are
\begin{equation}
\{H_{1},H_{3}\}=0, \qquad \{H_{2},H_{3}\}=H_{3}, \qquad
\{H_{1},H_{2}\}=2H_{1},
\end{equation}
hence it clearly solvable.  Taking a regular element $h\in {\mathcal G}^{\ast }$,
such that $h(H_{j})=h_{j}=0$,  for $j=1$ and $3$, and
$h(H_{2})=h_{2}\in  {\mathbb R}$ being arbitrary, one obtains the integrability of the
problem above in quadratures.

{\bf 0.2.} In 1974 V.~Arnold proved~[4] the following important result known as
the commutative (abelian) Liouville-Arnold theorem.

\medskip

\noindent
{\bf Theorem 0.5  (J.~Liouville -- V.~Arnold).} {\it Suppose a set
${\mathcal G}$ of functions $H_{j}:M^{2n}\rightarrow {\mathbb R}$,
$j=\overline{1,n}$, on a symplectic manifold $M^{2n}$  is abelian, that is
\begin{equation}
\{H_{i},H_{j}\}=0
\end{equation}
for all $i,j=\overline{1,n}$. If on the compact and connected
integral submanifold $M_{h}^{n} =\left\{u\in M^{2n} :\right.$
$\left. h(H_{j})=h_{j}\in {\mathbb R}, \ j=\overline{1,n}, \
h\in {\cal G}^{\ast }\right\}$ with $h\in {\mathcal G}$
being regular,  all functions $H:M^{2n}\in {\mathbb R}$, $j=\overline{1,n}$,
are functionally independent, then $M_{h}^{n}$ is diffeomorphic to the
$n$-dimensional torus ${\mathbb T}^{n}\simeq M^{2n}$, and the
motion on it with respect to the Hamiltonian $H=H_{1}\in {\mathcal G}$
is a quasi-periodic function of the evolution parameter $t \in {\mathbb R}$.}

\medskip

A dynamical system satisfying the hypotheses of Theorem 0.5 is called
completely integrable.

In 1978 Mishchenko and Fomenko [2] proved the following generalization of
the Liou\-vil\-le-Ar\-nold Theorem 0.5:

\medskip

\noindent
{\bf Theorem 0.6 (A.~Mishchenko -- A.~Fomenko).} {\it Assume that on a
symplectic manifold $\left(M^{2n},\omega^{(2)}\right)$ there is a
nonabelian Lie algebra ${\mathcal G}$ of invariants $H_{j}:M\in
{\mathbb R}$, $j=\overline{1,k}$, with respect to the dual
Poisson bracket on $M^{2n}$, that is
\begin{equation}
\{H_{i},H_{j}\}=\sum_{s=1}^{k}c_{ij}^{s}H_{s},
\end{equation}
where all values $c_{ij}^{s}\in {\mathbb R}$, $i,j,s =\overline{1,k}$,
are constants, and the following conditions are satisfied:

i) the integral submanifold $M_{h}^{r}
:=\left\{u\in M^{2n} :h(H_{j})=h_{j}\in {\mathbb R}, \ j=\overline{1,k},
\ h\in {\mathcal G}^{\ast }\right\}$ is compact and connected at a regular element
$h\in {\mathcal G}^{\ast }$;

ii) all functions $H_{j} :M^{2n}\rightarrow {\mathbb R}$, $j=\overline{1,k}$,
are functionally independent on $M^{2n}$;

iii) the Lie algebra ${\mathcal G}$ of invariants satisfies the
following relationship:
\begin{equation}
\dim {\mathcal G}+\mbox{\rm rank}\; {\mathcal G}=\dim M^{2n},
\end{equation}
where $\mbox{\rm rank}\, {\mathcal G}=\dim {\mathcal G}_{h}$ is the dimension
of a Cartan subalgebra ${\mathcal G}_{h}\subset {\mathcal G}$. Then the
submanifold $M_{h}^{r}\subset M^{2n}$ is
$r=\mbox{\rm rank}\; {\mathcal G}$-dimensional, invariant with respect each vector field
$K\in \Gamma (M^{2n})$, generated by an element
$H\in {\mathcal G}_{h}$, and diffeomorphic to the $r$-dimensional torus
${\mathbb T}^{r}\simeq M_{h}^{r}$,
on which the motion is a quasiperiodic function of the evolution
parameter $t\in {\mathbb R}$.}

\medskip

{\bf 0.3.} The simplest proof of the Mishchenko-Fomenko Theorem 0.6 can be
obtained from the well known [3, 16] classical Lie-Cartan theorem.

\medskip

\noindent
{\bf Theorem 0.7 (S.~Lie -- E.~Cartan).} {\it Suppose that a point $h\in
{\mathcal G}^{\ast }$ for a given Lie algebra ${\mathcal G}$
of invariants $H_{j}: M^{2n}\rightarrow {\mathbb R}$, $j=\overline{1,k}$,
is not critical, and the $\mbox{\rm rank}\,
||\{H_{i},H_{j}\}: i,j=\overline{1,k}\,||=
2(n-r)$ is constant in an open neighborhood $U_{h} \in
{\mathbb R}^{n}$  of the point $\left\{h(H_{j})=h_{j}\in {\mathbb R}:
j=\overline{1,k}\,\right\}\subset {\mathbb R}^{k}$.
Then in the neighborhood $(h\circ  H^{-1}:U_{h}
\subset M^{2n})$ there exist $k \in {\mathbb Z}_{+}$ independent functions
$f_{s}:{\mathcal G}\rightarrow {\mathbb R}$, $s=\overline{1,k}$,
such that the functions $F_{s}:=(f_{s}\circ H):M^{2n}\in {\mathbb R}$,
$s=\overline{1,k}$, satisfy the following relationships:
\begin{equation}
\{F_{1},F_{2}\}=\{F_{3},F_{4}\}= \cdots =\{F_{2(n-r)-1},F_{2(n-r)}\}=1,
\end{equation}
with  all  other brackets $\{F_{i},F_{j}\}=0$, where $(i,j)\neq (2s-1,2s)$,
$s=\overline{1,n-r}$.  In particular, $(k+r-n)\in {\mathbb Z}_{+}$
functions $F_{j}:M^{2n}\rightarrow {\mathbb R}$, $j=\overline{1,n-r}$,
and $F_{s}:M^{2n}\rightarrow {\mathbb R}$, $s=\overline{1,k-2(n-r)}$,
compose an abelian algebra ${\mathcal G}_{\tau }$
of new invariants on $M^{2n}$, independent on
$(h\circ H)^{-1}(U_{h})\subset M^{2n}.$}

\medskip

As a simple corollary of the Lie-Cartan Theorem 0.7 one obtains the
following: in the case of the Mishchenko-Fomenko theorem when
$\mbox{rank}\; {\mathcal G}+\dim{\mathcal G}=\dim M^{2n}$, that is $r+k=2n$,
the abelian algebra ${\mathcal G}_{\tau }$ (it is not a subalgebra
of ${\mathcal G}$!) of invariants on~M$^{2n}$
is~just $n=1/2\dim M^{2n}$-dimensional, giving rise to its local complete
integrability in \mbox{$(h\circ H)^{-1}(U_{h})\subset  M^{2n}$} via the abelian
Liouville-Arnold Theorem~0.5. It is also evident that the Mishchenko-Fomenko
nonabelian integrability Theorem~0.6 reduces to the commutative (abelian)
Liouville-Arnold case when a Lie algebra ${\mathcal G}$ of invariants is just
abelian, since then $\mbox{rank}\; {\mathcal G}=\dim{\mathcal G}=
1/2\dim M^{2n}=n\in {\bf Z}_{+} $ -- the standard  complete integrability condition.

All the cases of integrability by quadratures described above  pose the
following fundamental question: How can one ef\/fectively construct by means
of algebraic-analytical methods the corresponding integral submanifold
imbedding
\begin{equation}
{\pi }_{h}:M_{h}^{r}\rightarrow M^{2n},
\end{equation}
where $r=\dim \mbox{rank} \;{\mathcal G}$, thereby making it possible to express the
solutions of an integrable f\/low on $M_{h}^{r}$ as some exact quasi-periodic
functions on the torus ${\mathbb T}^{r}\simeq M_{h}^{r}$.

Below we shall describe an algebraic-analytical algorithm for resolving this
question for the case when a symplectic manifold $M^{2n}$ is
dif\/feomorphic to the canonically symplectic cotangent phase space
$T^{\ast }({\mathbb R})\simeq M^{2n}$.

\section{General setting}

\setcounter{equation}{0}

{\bf 1.1.} Our main object of study will be dif\/ferential
systems of vector f\/ields on the cotangent phase space $M^{2n}=T^{\ast }
({\mathbb R}^{n})$, $n\in {\mathbb Z}_{+}$, endowed with the canonical
symplectic structure $\omega ^{(2)}\in \Lambda ^{2}(M^{2n})$,  where by
$\omega ^{(2)}=d(\mbox{pr}^{\ast }\alpha^{(1)})$, and
\begin{equation}
\alpha ^{(1)}:=\langle p,dq\rangle =\sum_{j=1}^{n}p_{j}dq_{j},
\end{equation}
is the canonical 1-form on the base space
${\mathbb R}^{n}$, lifted naturally to the space
$\Lambda^{1}(M^{2n})$, $(q,p)\in M^{2n}$
are canonical coordinates on $T^{\ast }({\mathbb R}^{n})$,
$\mbox{pr}:T^{\ast }({\mathbb R}^{n})\rightarrow {\mathbb R}$  is
the canonical projection, and  $\langle \cdot ,\cdot \rangle $ is the usual scalar
product in ${\mathbb R}^{n}$.

Assume further that there is also given  a Lie subgroup $G$  (not
necessarily compact), acting symplectically via the mapping $\varphi
:G\times M^{2n}\rightarrow M^{2n}$ on $M^{2n},$ generating  a Lie algebra
homomorphism  $\varphi _{\ast }:T({\mathcal G)\rightarrow }\Gamma (M^{2n})$ via
the diagram
\begin{equation}
\begin{array}{ccccc}
{\mathcal G} & \simeq & T({\mathcal G}) & \stackrel{\varphi _{\ast }(u)}{\rightarrow
} & T(M^{2n}) \\
&  & \left\downarrow {}\right. &  & \left\downarrow {}\right. \\
&  & {\mathcal G} & \stackrel{\varphi (u)}{\rightarrow } & M^{2n}
\end{array}
\end{equation}
where $u\in M^{2n}.$ Thus, for any  $a\in {\mathcal G}$ one can
def\/ine a vector f\/ield $K_{a}\in \Gamma (M^{2n})$ as follows:
\begin{equation}
K_{a}=\varphi _{\ast }\cdot a.
\end{equation}
Since the manifold $\ M^{2n}$  is symplectic, one can naturally def\/ine
for any $a\in {\mathcal G}$ a function $H_{a}\in {\mathcal D(}M^{2n})$
as follows:
\begin{equation}
-i_{K_{a}}\omega ^{(2)}=dH_{a},
\end{equation}
whose existence follows from the invariance property
\begin{equation}
L_{K_{a}}\omega ^{(2)}=0
\end{equation}
for all $a\in {\mathcal G}.$  The following lemma [1] is useful in
applications.

\medskip

\noindent
{\bf Lemma 1.1.} {\it If the first homology group $H_{1}({\mathcal G};
{\mathbb R})$ of the Lie algebra ${\mathcal G}$ vanishes, then the mapping
$\Phi :{\mathcal G} \rightarrow {\mathcal D}(M^{2n})$ defined as
\begin{equation}
\Phi (a):=H_{a}
\end{equation}
for any $a\in {\mathcal G}$,  is  a Lie algebra homomorphism of
${\mathcal G}$ and  ${\mathcal D}(M^{2n})$  (endowed with the
Lie structure induced  by the symplectic structure
$\omega^{(2)}\in \Lambda ^{2}(M^{2n})$). In this case ${\mathcal G}$ is said to be
Poissonian.}

\medskip

As the mapping $\Phi :{\mathcal G} \rightarrow {\mathcal D}(M^{2n})$  is
evidently linear in ${\mathcal G}$, the expression (1.6) naturally def\/ines a
momentum mapping  $l:M^{2n}\rightarrow {\mathcal G}^{\ast }$  as follows: for
any $u\in M^{2n}$ and all  $a\in {\mathcal G}$
\begin{equation}
(l(u),a)_{{\mathcal G}}:=H_{a}(u),
\end{equation}
where $(\cdot ,\cdot )_{{\mathcal G}}$ is the standard scalar product on the
dual pair ${\mathcal G}^{\ast }\times {\mathcal G}.$ The following characteristic
equivariance [1] lemma holds.

\medskip

\noindent
{\bf Lemma 1.2.} {\it The diagram
\begin{equation}
\begin{array}{ccc}
M^{2n} & \stackrel{l}{\rightarrow } & {\mathcal G}^{\ast } \\
 \varphi _{g}\Bigr\downarrow &  & \Bigl\downarrow \mbox{\rm Ad}_{g^{-1}}^{\ast
} \\
M^{2n} & \stackrel{l}{\rightarrow } & {\mathcal G}^{\ast }
\end{array}
\end{equation}
commutes for all $g\in G$, where $\mbox{\rm Ad}_{g^{-1}}^{\ast }:
{\mathcal G}^{\ast }\rightarrow {\mathcal G}^{\ast }$  is the corresponding
co-adjoint action of the Lie group $G$ on the dual space ${\mathcal G}^{\ast }.$}

\medskip

Take now any vector $h\in {\mathcal G}^{\ast }$ and consider a subspace
${\mathcal G}_{h}\subset {\mathcal G}$, consisting of elements
$a\in {\mathcal G}$, such that $\mbox{ad}_{a}^{\ast }h=0,$  where
$\mbox{ad}_{a}^{\ast }:{\mathcal G}^{\ast}\rightarrow {\mathcal G}^{\ast }$
is the corresponding  Lie algebra  ${\mathcal G}$
representation in the dual space ${\mathcal G}^{\ast }$.

The following lemmas hold.

\medskip

\noindent
{\bf Lemma 1.3.} {\it The subspace ${\mathcal G}_{h}\subset {\mathcal G}$
is a Lie subalgebra of ${\mathcal G}$,  called here a Cartan subal\-geb\-ra.}

\medskip

\noindent
{\bf Lemma 1.4.} {\it Assume a vector $h\in {\mathcal G}^{\ast }$
is chosen in such a way that $r=\dim{\mathcal G}_{h}$ is minimal. Then the
Cartan Lie subalgebra ${\mathcal G}_{h}\subset {\mathcal G}$ is abelian.}

\medskip

In Lemma 1.4 the corresponding element  $h\in {\mathcal G}^{\ast }$ is called
regular and the number $r=\dim $ ${\mathcal G}_{h}$ is called  the
$\mbox{rank}\; {\mathcal G}$  of the Lie algebra ${\mathcal G}.$

{\bf 1.2.} Some twenty years ago Mishchenko and Fomenko~[2] proved the following
important noncommutative (nonabelian) Liouville-Arnold theorem.

\medskip

\noindent
{\bf Theorem 1.5.} {\it On a symplectic space $\left(M^{2n},
\omega ^{(2)}\right)$
let there be given a set of smooth func\-tions $H_{j}\in {\mathcal D}
(M^{2n})$, $j=\overline{1,k}$, whose linear span over ${\mathbb R}$
comprises a Lie algebra ${\mathcal G}$ with respect to the
corresponding Poisson bracket on $M^{2n}$.  Suppose also that the set
\[
M_{h}^{2n-k}:=\left\{u\in M^{2n}:h(H_{j})=h_{j}\in {\mathbb R}, \ j=\overline{1,k},
\ h\in {\mathcal G}^{\ast }\right\}
\]
with $h\in {\mathcal G}^{\ast }$ regular, is a submanifold of $M^{2n}$,
and on $M_{h}^{2n-k}$ all the functions $H_{j}\in {\mathcal D}(M^{2n})$,
$ j=\overline{1,k}$, are functionally independent.
Assume also that the Lie algebra ${\mathcal G}$  satisfies the following
condition:
\begin{equation}
\dim {\mathcal G}+\mbox{\rm rank}\; {\cal G} =\dim M^{2n}.
\end{equation}
Then the submanifold $M_{h}^{r}:=M_{h}^{2n-k}$ is $\mbox{\rm rank}\;
{\mathcal G}=r$-dimensional and invariant with respect to each vector field
 $K_{\overline{a}}\in \Gamma (M^{2n})$ with $\overline{a}\in {\mathcal G}_{h}
\subset {\mathcal G}$. Given a vector field $K=K_{\overline{a}}
\in \Gamma (M^{2n})$  with $\overline{a}\in {\mathcal G}_{h}$  or
$K\in \Gamma (M^{2n})$ such that  $[K,K_{a}]=0$ for all $a\in
{\mathcal G}$, then, if the submanifold $M_{h}^{r}$ is
connected and compact, it is diffeomorphic to the $r$-dimensional
torus ${\mathbb T}^{r}\simeq M_{h}^{r}$ and  the motion of the
vector field $K\in \Gamma (M^{2n})$ on it is a quasiperiodic function
of the evolution parameter $t\in {\mathbb R}$.}

\medskip

The easiest proof of this result can be obtained from the well known~[3]
classical  Lie-Cartan theorem, mentioned in the Introduction. Below  we
shall only sketch the original Mishchenko-Fomenko proof which is heavily
based on symplectic theory techniques, some of which have been discussed
above.

$\blacktriangleleft$ {\bfseries \itshape
Sketch of the proof.} Def\/ine a Lie group $G$
naturally as $G=\exp {\mathcal G}$, where ${\mathcal G}$ is the Lie algebra of
functions $H_{j}\in {\mathcal D}(M^{2n})$, $j=\overline{1,k}$, in the theorem,
with respect to the Poisson bracket  $\{\cdot ,\cdot \}$ on $M^{2n}$.
Then  for an element $h\in {\mathcal G}^{\ast }$ and any
$a=\sum\limits_{j=1}^{k}c_{j}H_{j}\in {\mathcal G}$,
where $c_{j}\in {\mathbb R}$, $ j=\overline{1,k}$, the following equality
\begin{equation}
(h,a)_{{\mathcal G}}:=\sum_{j=1}^{k}c_{j}h(H_{j})=\sum_{j=1}^{k}c_{j}h_{j}
\end{equation}
holds.  Since all functions $H_{j}\in {\mathcal D}(M^{2n})$,  $j=\overline{1,k}$,
 are independent on the level submanifold $M_{h}^{r}\subset M^{2n}$, this
evidently means that the element $h\in {\mathcal G}^{\ast }$ is regular  for
the Lie algebra  ${\mathcal G}$. Consequently, the Cartan Lie subalgebra
${\mathcal G}_{h}$ $\subset $ ${\mathcal G}$
is abelian. The latter is proved by
means of simple straightforward calculations. Moreover, the corresponding
momentum mapping $l:M^{2n}\rightarrow {\mathcal G}^{\ast }$  is constant on
$M_{h}^{r}$ and satisf\/ies the following relation:
\begin{equation}
l(M_{h}^{r})=h\in {\mathcal G}^{\ast }.
\end{equation}

\noindent
From this it can be shown that all vector f\/ields $K_{\overline{a}}\in
\Gamma (M^{2n})$, $\overline{a}\in {\mathcal G}_{h}$, are tangent to the
submanifold $M_{h}^{r}\subset M^{2n}$. Thus the corresponding Lie subgroup
$G_{h}:=\exp {\mathcal G}_{h}$ acts naturally and invariantly on $M_{h}^{r}$.
If the submanifold $M_{h}^{r}\subset M^{2n}$ is connected and compact, it
follows from (1.9) that  $\dim M_{h}^{r}=\dim M^{2n}-\dim {\mathcal G}=
\mbox{rank}\; {\mathcal G}=r$, and one obtains via the Arnold theorem~[4],
that $M_{h}^{r}\simeq {\mathbb T}^{r}$  and the motion of the vector f\/ield $K\in
\Gamma (M^{2n})$ is a quasiperiodic function of the evolution parameter
$t\in {\mathbb R}$; thus proving the theorem.~$\blacktriangleright $

As a nontrivial consequence of the Lie-Cartan theorem mentioned before and
of the Theorem~1.5, one can prove the following dual theorem about abelian
Liouville-Arnold integrability.

\medskip

\noindent
{\bf Theorem 1.6.} {\it Let a vector field $K\in \Gamma (M^{2n})$
be completely integrable via the nonabelian scheme of Theorem~1.5. Then it
is also Liouville-Arnold integrable on $M^{2n}$ and possesses,
under some additional conditions,  yet another abelian Lie algebra
${\mathcal G}_{h}$ of  functionally independent invariants on $M^{2n}$,
 for which $\dim {\mathcal G}_{h}=n=1/2\dim M^{2n}.$}

\medskip

The available proof of the theorem above is quite complicated, and we shall
comment on it in detail later on. We mention here only that some analogs of
the reduction Theorem~1.5  for the case where $M^{2n}\simeq
{\mathcal G}^{\ast }$, so that an arbitrary Lie group
$G$ acts symplectically on the manifold, were proved also
in [6--10, 34]. Notice  here, that in case
when the equality (1.9) is not satisf\/ied, one can then construct in the
usual way the reduced manifold  $\overline{M}_{h}^{2n-k-r}:=M_{h}^{2n-k}/G_{h}$
on which there exists a symplectic
structure $\overline{\omega }_{h}^{(2)}\in \Lambda ^{2}
\left(\overline{M}_{h}^{2n-k-r}\right)$, def\/ined  as
\begin{equation}
r_{h}^{\ast }\overline{\omega }_{h}^{(2)}=\pi _{h}^{\ast }\omega ^{(2)}
\end{equation}
with respect to the following compatible reduction-imbedding diagram:
\begin{equation}
\overline{M}_{h}^{2n-k-r}\stackrel{r_{h}}{\longleftarrow }M_{h}^{2n-k}
\stackrel{\pi _{h}}{\rightarrow }M^{2n},
\end{equation}
where  $r_{h}:M_{h}^{2n-k}\rightarrow \overline{M}_{h}^{2n-k-r}$  and
 $\pi _{h}:M_{h}^{2n-k}\rightarrow M^{2n}$ are, respectively, the
corresponding reductions and imbedding mappings. The nondegeneracy of the
2-form $\overline{\omega }_{h}^{(2)}\in \Lambda ^{2}(\overline{M}_{h})$
def\/ined by (1.12), follows simply from  the expression
\begin{equation}
\ba{l}
\ds \ker \left(\pi _{h}^{\ast }\omega ^{(2)}(u)\right)=
T_{u}\left(M_{h}^{2n-k}\right)\cap T_{u}^{\perp}\left(M_{h}^{2n-k}\right)
\vspace{3mm}\\
\ds \qquad =\mbox{span}_{{\mathbb R}}\left\{K_{\overline{a}}(u)
\in T_{u}\left(\overline{M}_{h}^{2n-k-r}:=M_{h}^{2n-k}/G_{h}\right):
\overline{a}\in {\mathcal G}_{h}\right\}
\ea
\ee
for any $u\in M_{h}^{2n-k}$,
since all vector f\/ields $K_{\overline{a}}\in \Gamma (M^{2n})$,
 $\overline{a}\in {\mathcal G}_{h}$, are tangent to
$\overline{M}_{h}^{2n-k-r}:=M_{h}^{2n-k}/G_{h}$.
Thus, the reduced space
$\overline{M}_{h}^{2n-k-r}:=M_{h}^{2n-k}/G_{h}$  with respect to the orbits
of the Lie subgroup $G_{h}$  action on  $M_{h}^{2n-k}$ will be a $(2n-k-r)$
-dimensional symplectic manifold. The latter evidently means that the number
$2n-k-r=2s\in {\mathbb Z}_{+}$  is even as  there is no symplectic structure
on odd-dimensional manifolds. This obviously is closely connected with the
problem of existence of a symplectic group action of a Lie group $G$ on a
given symplectic manifold $\left(M^{2n},\omega^{(2)}\right)$
with a symplectic structure $\omega ^{(2)}\in \Lambda ^{(2)}(M^{2n})$ being
apriori f\/ixed.  From this point of view one can consider the inverse
problem of constructing symplectic structures on a manifold  $M^{2n}$,
admitting a Lie group $G$ action. Namely, owing to the equivariance property
(1.8) of the momentum mapping $l:M^{2n}\rightarrow {\mathcal G}^{\ast }$, one
can obtain the induced symplectic structure $l^{\ast }\Omega _{h}^{(2)}\in
\Lambda ^{2}\left(\overline{M}_{h}^{2n-k-r}\right)$  on  $\overline{M}_{h}^{2n-k-r}$
from \ the canonical symplectic structure $\Omega _{h}^{(2)}\in \Lambda
^{(2)}(Or(h;G))$ on the orbit  $Or(h;G)\subset {\cal G}^{\ast }$ of a
regular element $h\in {\mathcal G}^{\ast }$.
Since the symplectic structure
$l^{\ast }\Omega _{h}^{(2)}\in \Lambda ^{2}(\overline{M}_{h})$ can be
naturally lifted to the 2-form  $\widetilde{\omega }^{(2)}=(r_{h}^{\ast
}\circ l^{\ast })\Omega _{h}^{(2)}\in \Lambda ^{2}\left(M_{h}^{2n-k}\right)$,
the latter being degenerate on  $M_{h}^{2n-k}$
can apparently be nonuniquely extended
on the whole manifold $M^{2n}$  to a symplectic structure  $\omega
^{(2)}\in \Lambda ^{2}(M^{2n}),$ for which the action of the Lie group
$G$ is  apriori  symplectic. Thus, many properties of a
given dynamical system with  a Lie algebra  ${\mathcal G}$ of invariants on
$M^{2n}$ are deeply connected with the symplectic structure $\omega ^{(2)}\in
\Lambda ^{2}(M^{2n})$ the  manifold $M^{2n}$ is endowed with, and in
particular, with the corresponding integral submanifold imbedding mapping
$\pi _{h}:M_{h}^{2n-k}\rightarrow M^{2n}$ at a regular element $h\in
{\mathcal G}^{\ast }$.
The problem of direct algebraic-analytical construction
of this mapping was in part solved in~[11] in the case where $n=2$  for
an abelian algebra ${\mathcal G}$ on the manifold
$M^{4}=T^{\ast }({\mathbb R}^{2})$.
The treatment  of this problem in~[11] has been
extensively based both on the classical Cartan studies of integral
submanifolds of ideals in Grassmann algebras and on the modern
Galisau-Reeb-Francoise results for a symplectic manifold
$\left(M^{2n},\omega^{(2)}\right)$ structure, on which there
 exists an involutive set  ${\mathcal G}$
of functionally independent invariants  $H_{j}\in {\mathcal D}(M^{2n})$,
\mbox{$ j=\overline{1,n}$}. In what follows below we generalize the
Galisau-Reeb-Francoise  results to the case of a nonabelian set of
functionally independent functions  $H_{j}\in {\mathcal D}(M^{2n})$,
\mbox{$j=\overline{1,k}$}, comprising a Lie algebra  ${\mathcal G}$  and
satisfying the Mishchenko-Fomenko condition~(1.9):
$\dim {\mathcal G}+\mbox{rank} \; {\mathcal G} =\dim M^{2n}$.
This makes it possible to devise an ef\/fective
algebraic-analytical method of constructing  the corresponding integral
submanifold imbedding and reduction mappings, giving rise to a wide class of
exact, integrable by quadratures solutions of a given integrable vector
f\/ield on $M^{2n}$.

\section{Integral submanifold imbedding problem for an
abelian Lie algebra of invariants}

\setcounter{equation}{0}

{\bf 2.1.} We shall consider  here only a set ${\mathcal G}$
 of commuting polynomial functions $H_{j}\in {\mathcal D}(M^{2n})$,
$j=\overline{1,n}$, on the canonically symplectic phase space
 $M^{2n}=T^{\ast}({\mathbb R}^{n})$.
Due to the Liouville-Arnold theorem~[4], any
dynamical system $K\in \Gamma (M^{2n})$ commuting  with corresponding
Hamiltonian vector f\/ields $K_{a}$ for all $a\in {\mathcal G}$,
will be integrable by quadratures in case of a regular element
 $h\in {\mathcal G}^{\ast}$, which def\/ines the corresponding
integral submanifold $M_{h}^{n}:=\left\{u\in M^{2n}:h(H_{j})=h_{j}
\in {\mathbb R}, \ j=\overline{1,n}\,\right\}$
which is dif\/feomorphic (when compact and connected) to
the $n$-dimensional torus ${\mathbb T}^{n}\simeq M_{h}^{n}$.
This in particular means that there exists some algebraic-analytical
expression for the integral submanifold imbedding mapping
$\pi _{h}:M_{h}^{n}\rightarrow M^{2n}$
into the ambient phase space $M^{2n}$, which one should f\/ind
in order to properly demonstrate integrability by quadratures.

The problem formulated above was posed and in part solved (as was mentioned
above) for $n=2$ in~[11] and in [13] for a Henon-Heiles dynamical
system which had previously been integrated [14, 15] using other tools. Her
we generalize the approach of [11] for the general case
$n\in {\mathbb Z}_{+}$ and proceed further in Chapter~3 to solve
this problem in the case of
a nonabelian Lie algebra ${\mathcal G}$  of polynomial invariants on
$M^{2n}=T^{\ast }({\mathbb R}^{n})$, satisfying all the conditions of
Mishchenko-Fomenko Theorem~1.5.

{\bf 2.2.} Def\/ine now the basic vector f\/ields $K_{j}\in \Gamma (M^{2n})$,
$j=\overline{1,n}$, generated by basic elements $H_{j}\in {\mathcal G}$
of  an abelian Lie algebra ${\mathcal G}$ of invariants on $M^{2n}$, as
follows:
\begin{equation}
-i_{K_{j}}\omega ^{(2)}=dH_{j}
\end{equation}
for all $j=\overline{1,n}$. It is easy to see that the condition
$\{H_{j},H_{i}\}=0$ for all
$i,j=\overline{1,n\ }$,  yields also
$[K_{i},Kj]=0$  for all
$i,j=\overline{1,n}.$ Taking into account that
$\dim M^{2n}=2n$ one obtains the equality
$(\omega ^{(2)})^{n}=0$ identically on
$M^{2n}$. This makes it possible  to formulate the
following Galisau-Reeb result.

\medskip

\noindent
{\bf Theorem 2.1.} {\it Assume that an element $h\in {\mathcal G}^{\ast }$
is chosen to be regular and a Lie algebra ${\mathcal G}$
of invariants on $M^{2n}$ is abelian. Then there exist
differential 1-forms $h_{j}^{(1)}\in \Lambda ^{1}(U(M_{h}^{n}))$,
$j=\overline{1,n}$, where $U(M_{h}^{n})$ is some
open neighborhood of the integral submanifold $M_{h}^{n}\subset M^{2n}$,
satisfying  the following properties:

i) $\left. \omega^{(2)}\right|_{U(M_{h}^{n})}=
\sum\limits_{j=1}^{n}dH_{j}\wedge h_{j}^{(1)};$

ii) the exterior differentials $dh_{j}^{(1)}\in \Lambda
^{2}(U(M_{h}^{n}))$ belong to the ideal ${\mathcal I}({\mathcal G})$
 in the Grassmann algebra $\Lambda (U(M_{h}^{n}))$,
generated  by $1$-forms $dH_{j}\in \Lambda^{1}(U(M_{h}^{n}))$,
$j=\overline{1,n}$.}

\medskip

$\blacktriangleleft$ {\bfseries \itshape Proof.}
Consider the following identity on $M^{2n}$:
\begin{equation}
(\otimes _{j=1}^{n}i_{K_{j}})\left(\omega ^{(2)}\right)^{n+1}=
0=\pm (n+1)!\left(\wedge _{j=1}^{n}dH_{j}\right)\wedge  \omega ^{(2)},
\end{equation}
 which implies that the 2-form  $\omega ^{(2)}\!\in\! {\mathcal I}
({\mathcal G})$. Whence, one can f\/ind 1-forms $h_{j}^{(1)}\!\in\!
\Lambda ^{1}(U(M_{h}^{n}))$, $j=\overline{1,n}$,
satisfying the condition
\begin{equation}
\left. \omega ^{(2)}\right| _{U(M_{h}^{n})}=\sum_{j=1}^{n}dH_{j}\wedge
h_{j}^{(1)}.
\end{equation}
Since $\omega ^{(2)}\in \Lambda ^{2}(U(M_{h}^{n}))$ is
nondegenerate on $M^{2n}$, it follows that all 1-forms
$h_{j}^{(1)}$, $j=\overline{1,n}$, in (2.3) are independent on
$U(M_{h}^{n})$, proving part  $i)$  of the theorem.
As $d\omega ^{(2)}=0$ on $M^{2n}$, from (2.3) one gets that
\begin{equation}
\sum_{j=1}^{n}dH_{j}\wedge dh_{j}^{(1)}=0
\end{equation}
on $U(M_{h}^{n})$, hence it is obvious that $dh_{j}^{(1)}\in
{\mathcal I}({\mathcal G })\subset \Lambda (U(M_{h}^{n}))$  for all
$j=\overline{1,n}$, proving part  $ii)$ of the theorem.~$\blacktriangleright$

Now we proceed to study properties of the integral submanifold
$M_{h}^{n}\subset M^{2n}$ of the ideal ${\mathcal I}({\mathcal G})$ in the
Grassmann algebra  $\Lambda (U(M_{h}^{n}))$. In general, the integral
submanifold $M_{h}^{n}$ is completely described~[16] by means of the
imbedding
\begin{equation}
\pi _{h}:M_{h}^{n}\rightarrow M^{2n}
\end{equation}
and using this, one can reduce all vector f\/ields $K_{j}\in \Gamma
(M^{2n})$, $j=\overline{1,n}$, on the submanifold
$M_{h}^{n}\subset M^{2n}$, since they are all evidently in its tangent space.
 If  $\overline{K}_{j}\in \Gamma (M_{h}^{n})$,
$j=\overline{1,n}$, are the corresponding pulled-back vector f\/ields
$K_{j}\in \Gamma (M^{2n})$,
$j=\overline{1,n}$, then by def\/inition, the equality
\begin{equation}
\pi _{h \ast }\circ \overline{K}_{j}=K_{j}\circ \pi _{h}
\end{equation}
holds for all $j=\overline{1,n}$. Similarly one can construct
1-forms $\overline{h}_{j}^{(1)}:=\pi _{h}^{\ast }\circ h_{j}^{(1)}\in
\Lambda ^{1}(M_{h}^{n})$, $j=\overline{1,n}$, which are
characterized by the following Cartan-Jost~[16] theorem.

\medskip

\noindent
{\bf Theorem 2.2.} {\it The following assertions are true:

i)  the 1-forms $\overline{h}_{j}^{(1)}\in \Lambda
^{1}(M_{h}^{n})$, $j=\overline{1,n}$, are independent on
$M_{h}^{n}$;

ii) the 1-forms $\overline{h}_{j}^{(1)}\in \Lambda ^{1}(M_{h}^{n})$,
 $j=\overline{1,n}$, are exact on $M_{h}^{n}$ and
satisfy $\overline{h}_{j}^{(1)}(\overline{K}_{j})=\delta _{ij}$,
$i,j=\overline{1,n}$.}

\medskip

$\blacktriangleleft$ {\bfseries \itshape Proof.} As the ideal
${\mathcal I}({\mathcal G})$ is by
def\/inition vanishing on $M_{h}^{n}\subset M^{2n}$ and closed on
$U(M_{h}^{n})$, the integral submanifold $M_{h}^{n}$ is well def\/ined
in the case of a regular element $h\in {\mathcal G}^{\ast }$. This implies
that the imbedding (2.5)  is nondegenerate on  $M_{h}^{n}\subset M^{2n}$,
or the 1-forms  $\overline{h}_{j}^{(1)}:=\pi _{h}^{\ast }\circ h_{j}^{(1)}$,
$j=\overline{1,n}$, will persist in being independent if they are
1-forms  $h_{j}^{(1)}\in \Lambda ^{1}(U(M_{h}^{n}))$,
 $j=\overline{1,n}$, proving part {\it i)}
 of the theorem. Using property {\it ii)} of
Theorem~2.1, one sees that on the integral submanifold
$M_{h}^{n}\subset M^{2n}$ all 2-forms $d\overline{h}_{j}^{(1)}=0$,
$j=\overline{1,n}$. Consequently, owing to the Poincar\'{e} lemma [1, 16],
the 1-forms  $\overline{h}_{j}^{(1)}=d\overline{t}_{j}\in \Lambda
^{1}(M_{h}^{n})$, $j=\overline{1,n}$, for some mappings
$\overline{t}_{j}:M_{h}^{n}\rightarrow {\mathbb R}$,
$j=\overline{1,n}$, def\/ining global coordinates
on an appropriate universal covering of $M_{h}^{n}$.
Consider now the following identity based on the representation~(2.3):
\begin{equation}
\left. i_{K_{j}}\omega ^{(2)}\right|
_{U(M_{h}^{n})}=-\sum_{i=1}^{n}h_{i}^{(1)}(K_{j}) dH_{i}:=-dH_{j},
\end{equation}
which holds for any $j=\overline{1,n}$. As all $dH_{j}\in \Lambda
^{1}(U(M_{h}^{n}))$, $j=\overline{1,n}$, are independent, from
(2.7) one infers that $h_{i}^{(1)}(K_{j})=\delta _{ij}$ for all
$i,j=\overline{1,n}$. Recalling now that for any $i=\overline{1,n}$,
$K_{i}\circ \pi _{h}=\pi _{h \ast }\circ K_{i}$,  one readily
computes that $\overline{h}_{i}^{(1)}(\overline{K}_{j})=\pi
_{h}^{\ast }h_{i}^{(1)}(\overline{K}_{j}):=h_{i}^{(1)}
\left(\pi _{h \ast }\circ K_{j}\right):=h_{i}^{(1)}(K_{j}\circ \pi _{h})=\delta _{ij}$  for all
$i,j=\overline{1,n}$,  proving part  {\it ii)} of the
theorem.~$\blacktriangleright$

The following is a simple consequence of Theorem 2.2:

\medskip

\noindent
{\bf Corollary 2.3.} {\it Suppose that the vector fields $K_{j}\in
\Gamma (M^{2n})$, $j=\overline{1,n}$, are parametrized globally along
their trajectories by means of the corresponding parameters
$t_{j}:M^{2n}\rightarrow {\mathbb R}$,
$j=\overline{1,n}$, that is on the phase space $M^{2n}$
\begin{equation}
d/dt_{j}:=K_{j}
\end{equation}
for all $j=\overline{1,n}$. Then the following important
equalities hold (up to constant normalizations) on the integral submanifold
$M_{h}^{n}\subset M^{2n}$:
\begin{equation}
\left. t_{j}\right| _{M_{h}^{n}}=\overline{t}_{j},
\end{equation}
where $1\leq j\leq n$.}

\medskip

{\bf 2.3.} To proceed with our investigation of the algebraic-analytical
properties of the imbedding (2.5), we shall summarize some useful facts
about the canonical transformations [1, 4] of symplectic manifolds and their
generating functions.

Suppose  that a symplectomorphism  $\Phi :M^{2n}\rightarrow \tilde{M}^{2n}$
satisf\/ies the condition
\begin{equation}
\Phi ^{\ast }\tilde{\omega}^{(2)}=\omega ^{(2)},
\end{equation}
where  $\tilde{\omega}^{(2)}\in \Lambda ^{2}(\tilde{M}^{2n})$ is a
symplectic structure on  $\tilde{M}^{2n}=T^{\ast}(\tilde{M}^{n})$. Since by
assumption $\omega ^{(2)}=d(\mbox{pr}^{\ast }\alpha ^{(1)})$,
where  $\alpha ^{(1)}\in \Lambda ^{1}({\mathbb R}^{n})$
is def\/ined by (1.1), and there exists a local 1-form
$\tilde{\alpha}^{(1)}\in \Lambda ^{1}(\tilde{M}^{n})$ such that
\begin{equation}
\mbox{pr}^{\ast }\alpha ^{(1)}-\mbox{pr}^{\ast }
\tilde{\alpha}^{(1)}+d<\widetilde{p},\, \widetilde{q}>=dS,
\end{equation}
where locally,  $S:{\mathbb R}^{n}\times {\mathbb R}^{n}\rightarrow
{\mathbb R}$ is a dif\/ferentiable
mapping called a generating function. Def\/ining on
$\tilde{M}^{n}$
\begin{equation}
\tilde{\alpha}^{(1)}=\sum_{j=1}^{n}\tilde{p}_{j}d\tilde{q}_{j},
\end{equation}
where $\tilde{p}\in T^{\ast }(\tilde{M}^{n})$ are canonical local (momentum)
 coordinates, from (2.11), (2.12) and (1.1) one readily f\/inds that
\begin{equation}
p_{j}=\partial S(q,\tilde{p})/\partial q_{j},\qquad
\tilde{q}_{j}=\partial S(q,\tilde{p})/\partial \tilde{p}_{j}
\end{equation}
for any  $j=\overline{1,n}$. The mapping $\Phi :M^{2n}\rightarrow
\tilde{M}^{2n}$ should in local coordinates satisfy the following condition:
\begin{equation}
\det (\partial \widetilde{p}(q,p)/\partial p)\neq 0
\end{equation}
almost everywhere (a.e.)  on $M^{2n}$. Since due to (2.14)  one can def\/ine
a.e. the mapping  $p:{\mathbb R}^{n}\times \tilde{M}^{n}\rightarrow
{\mathbb R}^{n}$,
the equations (2.13) determine the generating function $S:{\mathbb R}
^{n}\times {\mathbb R}^{n}\rightarrow {\mathbb R}$ up to some constant. And
conversely  [4], if there is given a generating function
 $S:{\mathbb R}^{n}\times {\mathbb R}^{n}\rightarrow {\mathbb R}$,
satisfying a.e. the condition
\begin{equation}
\det (\partial ^{2}S(q,\widetilde{p})/\partial q\partial p)\neq 0
\end{equation}
on ${\mathbb R}^{n}\times {\mathbb R}^{n}$, then one can determine a canonical
transformation $\Phi :M^{2n}\rightarrow \tilde{M}^{2n}$  of symplectic
manifolds $M^{2n}$ and $\tilde{M}^{2n}$, satisfying a.e. the condition~(2.14).

Assume now additionally that the submanifold $\tilde{M}^{n}\subset
\tilde{M}^{2n}$  coincides up to a dif\/feomor\-phism  with
the integral submanifold $M_{h}^{n}$  of the ideal  ${\mathcal I(G)}$
considered above. Then the corresponding symplectic manifold  $\tilde{M}
_{({\mathcal G}^{\ast })}^{2n}:=\cup _{h\in {\mathcal G}^{\ast }}T^{\ast
}(M_{h}^{n})$ is the usual topological sum of cotangent spaces, giving
rise to a natural f\/ibration of the symplectic manifold
$M^{2n}=T^{\ast }({\mathbb R}^{n})$:
\begin{equation}
M^{2n}=T^{\ast }({\mathbb R}^{n})\simeq \cup _{h\in {\mathcal G}^{\ast }}
T^{\ast }(M_{h}^{n})=\tilde{M}_{({\mathcal G}^{\ast })}^{2n}.
\end{equation}
The representation (2.16) appears to be very useful for treating the
imbedding (2.5). Namely, assume further that the integral submanifold
 $M_{h}^{n}\subset M^{2n}$ admits  a.e. on
$M_{h}^{n}$ coordinate charts from the base space
${\mathbb R}^{n}$ of the  entire phase space
$M^{2n}=T^{\ast }({\mathbb R}^{n})$. This
evidently means that the set  of annihilating 1-forms  $dH_{j}\in \Lambda
^{1}(M^{2n})$, $j=\overline{1,n}$, on $M_{h}^{n}$ must be solvable with
respect to the cotangent dif\/ferentials $dp_{j}\in \Lambda
^{1}(M_{h}^{n})$, $j=\overline{1,n}$:
\begin{equation}
\left\{\left. dH_{j}\right| _{M_{h}^{n}}=0,\
j=\overline{1,n}\, \right\}\Rightarrow \left\{dp_{j}=
\sum_{k=1}^{n}Q_{jk}(q,p)dq_{j}:(q,p)\in T^{\ast}(M_{h}^{n})\right\},
\end{equation}
where  $Q:  T^{\ast }(M_{h}^{n})\rightarrow \mbox{Hom}({\mathbb R}^{n})$ is
an invertible a.e. mapping. The implication (2.17) and a result of Arnold~[4]
imply the existence on $M_{h}^{n}$ of special cyclic coordinates
realizing  the isomorphisms  $M_{h}^{n}\simeq {\mathbb T}^{n}\simeq
\otimes _{j=1}^{n}{\mathbb S}_{j}^{1}$ (when
it is compact and connected as will be assumed in the sequel). It follows
directly from (2.16) that the phase space can be represented up to a
symplectomorphism as
\begin{equation}
M^{2n}=T^{\ast }({\mathbb R}^{n})\simeq \cup _{h\in {\mathcal G}^{\ast }}
T^{\ast}\left(\otimes _{j=1}^{n}{\mathbb S}_{j}^{1}\right)
\simeq \tilde{M}_{({\mathcal G}^{\ast })}^{2n}
\end{equation}
and the integral submanifold $M_{h}^{n}\subset M^{2n}$ can be covered by
charts with images in the base space  ${\mathbb R}^{n}\subset T^{\ast }
({\mathbb R}^{n})$ of the phase space $M^{2n}$. On the integral
submanifold $M_{h}^{n}$ there  is induced a canonical 1-form
$\alpha _{h}^{(1)}\in \Lambda ^{1}(M_{h}^{n})$ with respect to the
imbedding~(2.5), projected upon the base space ${\mathbb R}^{n}$:
\begin{equation}
\alpha _{h}^{(1)}:=\pi _{h}^{\ast }\circ \mbox{pr}^{\ast }\alpha ^{(1)}.
\end{equation}
This means that the following diagram is commutative:
\[
\begin{array}[b]{rrrrrrr}
T^{\ast }(T^{\ast }({\mathbb R}^{n})) &
\stackrel{\pi _{h}^{\ast }}{\rightarrow }
& T^{\ast }(M_{h}^{n}) & \left. \simeq \right| _{\mbox{\scriptsize loc}} & T^{\ast }
({\mathbb R}^{n}) &
\stackrel{\mbox{\scriptsize pr}^{\ast }}{\rightarrow } & T^{\ast }(T^{\ast }({\mathbb R}
^{n})) \\
\mbox{pr}^{\prime }\Bigr\downarrow  &  & \Bigl\downarrow \mbox{pr}_{h}
&  & \Bigl\downarrow \mbox{pr}  &  & \Bigl\downarrow \mbox{pr}^{\prime}
\\
M^{2n}=T^{\ast }({\mathbb R}^{n}) & \stackrel{\pi _{h}}{\longleftarrow } &
M_{h}^{n} & \left. \simeq \right| _{\mbox{\scriptsize loc}}
& {\mathbb R}^{n} & \stackrel{\mbox{\scriptsize pr}}{
\longleftarrow } & M^{2n}=T^{\ast }({\mathbb R}^{n})
\end{array},
\]
where $\mbox{pr}^{\prime}:T^{\ast }(T^{\ast }({\mathbb R}^{n}))\rightarrow
T^{\ast }({\mathbb R}^{n})$ is the standard projection mapping. The
representation (2.19) together with the isomorphism $M_{h}^{n}\simeq
\otimes _{j=1}^{n}{\mathbb S}_{j}^{1}$ implies, in particular, that in
circle-like coordinates $\mu \in \otimes _{j=1}^{n}{\mathbb S}_{j}^{1}$
the canonical 1-form $\alpha _{h}^{(1)}\in \Lambda ^{1}(U(M_{h}^{n}))$
can be written as
\begin{equation}
\alpha _{h}^{(1)}=\sum_{j=1}^{n}w_{j}d\mu _{j}=\pi _{h}^{\ast }\circ
\mbox{pr}^{\ast }\alpha ^{(1)},
\end{equation}
which is naturally lifted to a canonical 1-form $\mbox{pr}_{h}^{\ast }\circ
\alpha _{h}^{(1)}\in \Lambda ^{1}($\ $T^{\ast }(U(M_{h}^{n})))$,
where  $\mbox{pr}_{h}:T^{\ast }(M_{h}^{n})\rightarrow M_{h}^{n}$ is the natural
projection upon the base $M_{h}^{n}$, and
$w:\otimes _{j=1}^{n}{\mathbb S}_{j}^{1}\rightarrow {\mathbb R}^{n}$
is an smooth a.e. mapping parametrized by
$h\in {\mathcal G}^{\ast }$. Thus, in local coordinates on
$T^{\ast }(\otimes _{j=1}^{n}{\mathbb S}_{j}^{1})$
the imbedding (2.5) parametrized by regular elements
$h\in {\mathcal G}^{\ast }$ has the following form:
\begin{equation}
q_{j}=q_{j}(\mu ;h),\qquad p_{j}=p_{j}(\mu ;h),
\end{equation}
where in virtue of (2.20) for any $j=\overline{1,n}$
\begin{equation}
p_{j}= \sum_{i=1}^{n}w_{i}(\mu ;h)\partial \mu _{i}(q;h)/\partial
q_{j}\Bigr| _{M_{h}^{n}},
\end{equation}
with the mapping $q:=\mbox{pr}\circ \pi _{h}:M_{h}^{n}\rightarrow {\mathbb R}^{n}$
being invertible owing to the implication (2.17),
$\mu :{\mathbb R}^{n}\rightarrow \otimes _{j=1}^{n}{\mathbb S}_{j}^{1}$
is its inverse, and the mapping
$w:(\otimes _{j=1}^{n}{\mathbb S}_{j}^{1})\times {\mathcal G}^{\ast}
\rightarrow {\mathbb R}^{n}$ is as yet not def\/ined. The above analysis of
the imbedding problem for the integral submanifold  $M_{h}^{n}\subset
M^{2n}$ in the case when the implication (2.17) is solvable, tells us that
the Liouville-Arnold foliation (2.18) can be described ef\/fectively by
choosing a special parametrization by elements $h\in {\mathcal G}^{\ast }$,
for which the mapping  $w:(\otimes _{j=1}^{n}{\mathbb S}_{j}^{1})\times
{\mathcal G}^{\ast }\rightarrow {\mathbb R}^{n}$
is separable in the variables $\mu \in \otimes _{j=1}^{n}{\mathbb S}_{j}^{1}$,
that is on $M_{h}^{n}$ the expressions
\begin{equation}
w_{j}=w_{j}(\mu _{j};h)
\end{equation}
should hold for all $j=\overline{1,n}$. Such a case is called [1, 4]
the Hamilton-Jacobi separation of variables method which can now be
applied very naturally to our problem of f\/inding the embedding (2.5) by
algebraic-analytical means.

To proceed we f\/irst apply f\/irst to $\tilde{M}^{2n}$ a canonical
transformation into a new f\/ibration $\tilde{M}_{({\mathbb R}^{n})}^{2n}$
of the phase space $M^{2n}=T^{\ast }({\mathbb R}^{n})$ as def\/ined in (2.18),
def\/ined as follows:
\begin{equation}
M^{2n}=T^{\ast }({\mathbb R}^{n})\simeq \cup _{\tau \in {\mathbb R}^{n}}
T^{\ast}(M_{\tau }^{n}):=\tilde{M}_{({\mathbb R}^{n})}^{2n},
\end{equation}
where $M_{\tau }^{n}$ is an integral submanifold of a dual integrable
ideal ${\mathcal I}(h^{(1)})\subset \Lambda (U(M_{h}^{n}))/{\mathcal I}
({\mathcal G}_{\ }^{\ast })$, which is generated  by closed 1-forms $\tilde{h}
_{j}^{(1)}\in \Lambda ^{1}(U(M_{h}^{n}))/{\mathcal I}({\mathcal G}_{}^{\ast })$,
$j=\overline{1,n}$, obtained via an extension (preserving closedness) of
the closed 1-forms $\bar{h}_{j}^{(1)}\in \Lambda ^{1}(M_{h}^{n})$,
$j=\overline{1,n}$.
This obviously means that the constructed ideal ${\mathcal I}(h^{(1)})$
is also integrable on  $U(M_{h}^{n})\subset M^{2n}$, and
possesses integral submanifolds $M_{\tau }^{n}\subset U(M_{h}^{n})$
parametrized by a constant vector parameter $\tau:=(t_{1},t_{2},\ldots,
t_{n})\in {\mathbb R}^{n},$ composed of evolution parameters
 $t_{j}\in {\mathbb R}$, $j=\overline{1,n}$, of corresponding vector f\/ields
$K_{j}=d/dt_{j}$, $j=\overline{1,n}$, on the neighborhood
$U(M_{h}^{n}) \subset M^{2n}$. It follows that the two f\/ibrations
$\tilde{M}_{({\mathcal G}^{\ast })}^{2n}$ and
$\tilde{M}_{({\mathbb R}^{n})}^{2n}$  are locally
dif\/feomorphic in an open neighborhood $U(M_{h}^{n})$
of the integral submanifold  $M_{h}^{n}\subset M^{2n}$. As
they are symplectomorphic, one f\/inds that
\begin{equation}
\alpha _{\tau }^{(1)}=-\sum_{j=1}^{n}t_{j}dh_{j}=\pi _{\tau }^{\ast}
\mbox{pr}^{\ast }\alpha ^{(1)},
\end{equation}
holds on  the integral submanifold $M_{\tau }^{n}$ for the canonical
1-form $\alpha _{\tau }^{(1)}\in \Lambda ^{1}(M_{\tau }^{n})$, where
$\pi _{\tau }:M_{\tau }^{n}\rightarrow M^{2n}$ is the corresponding imbedding
mapping. Recall now that each symplectomorphism of the manifold  $M^{2n}$
is described [1, 4] by means of the corresponding generating function
$S:M_{h}^{n}\times {\mathbb R}^{n}\rightarrow {\mathbb R}$,
where by def\/inition,
\begin{equation}
\alpha _{h}^{(1)}=\alpha _{\tau }^{(1)}+dS.
\end{equation}
Making use now of the expressions (2.20) and (2.25), one obtains the local
relationship
\begin{equation}
\sum_{j=1}^{n}w_{j}d\mu _{j}+\sum_{j=1}^{n}t_{j}dh_{j}=dS(\mu ;h)
\end{equation}
on some open neighborhood  $M_{h,\tau }^{2n}\subset U(M_{h}^{n})\cap
U(M_{\tau }^{n})$,  where $U(M_{\tau }^{n})\subset M^{2n}$ is an open
neighborhood of the integral submanifold $M_{\tau }^{n}\subset M^{2n}$.

Using the Hamilton-Jacobi separability condition (2.23) in (2.27), one gets
\begin{equation}
S(\mu ;h)=\sum_{j=1}^{n}\int_{\mu _{j}^{0}}^{\mu _{j}}w_{j}(\lambda;h)
d\lambda
\end{equation}
for all $(\mu ;h)\in M_{h,\tau }^{2n}\subset M^{2n}$, where
$\mu ^{0}\in M_{h}^{n}\cap M_{h,\tau }^{2n}$ is some f\/ixed point. Whence one
obtains the desired expression for a vector parameter
$\tau :=(t_{1},t_{2},\ldots,t_{n})\in {\mathbb R}^{n}$:
\begin{equation}
t_{j}=\partial S(\mu ;h)/\partial h_{j},
\end{equation}
where  $j=\overline{1,n}$  and $(\mu ;h)\in M_{h,\tau }^{2n}$. On
the other hand, in virtue of (2.9) the set of parameters (2.29) can be
represented dually as dif\/ferential forms
\begin{equation}
\bar{h}_{j}^{(1)}=dt_{j}=\partial dS\left( \mu ;h\right) /\partial
h_{j}=\sum_{i=1}^{n}\left( \partial w_{i}(\mu _{i};h)/\partial h_{j}\right)
d\mu _{i}
\end{equation}
for all  $j=\overline{1,n}$  on $M_{h,\tau }^{2n}$. Since the 1-forms
$\bar{h}_{j}^{(1)}\in \Lambda ^{1}(M_{h,\tau }^{2n})$,
$j=\overline{1,n}$,
are assumed here to be known explicitly from the characteristic equation
(2.3), we can write them as
\begin{equation}
\bar{h}_{j}^{(1)}=\sum_{i=1}^{n}\bar{h}_{ji}(q,p)dq_{i},
\end{equation}
where $\bar{h}_{ji}:M_{h,\tau }^{2n}\rightarrow {\mathbb R}$,
$i,j=\overline{1,n}$ are some algebraic expressions.
Making use of the representations
(2.21) and (2.28) and equation (2.31), the set of 1-forms (2.30) is reduced
to the following purely dif\/ferential-algebraic relationships on
$M_{h,\tau }^{2n}$:
\begin{equation}
\partial w_{i}(\mu _{i};h)/\partial h_{j}={\bf P}_{ji}(\mu ,w;h),
\end{equation}
generalizing similar ones of [31, 18], where the characteristic functions
${\bf P}_{ji}:T^{\ast }(M_{h}^{n})\rightarrow {\mathbb R}$,
$i,j=\overline{1,n}$,  are def\/ined as follows:
\begin{equation}
{\bf P}_{ji}(\mu ,w;h):=\sum_{i=1}^{n}\bar{h}_{js}(q(\mu ;h),
w\partial \mu /\partial q(\mu ;h))\partial q_{s}/\partial \mu _{i}.
\end{equation}

A simple analysis of the relationships (2.32) and (2.33) tells us that for
all  $j=\overline{1,n}$  and $i\neq s=\overline{1,n}$ the following
algebraic relations
\begin{equation}
\partial {\bf P}_{ji}(\mu ;w;h)/\partial w_{s}=0
\end{equation}
must be satisf\/ied identically if $i\neq s$. It is clear that the above set
of purely dif\/ferential-algebraic relationships (2.33) and (2.34) makes it
possible to write down explicitly some f\/irst order compatible
dif\/ferential-algebraic equations, whose solution yields the f\/irst half of
the desired imbedding (2.5) for the integral submanifold $M_{h}^{n}\subset
M^{2n}$ in an open neighborhood $M_{h,\tau }^{2n}\subset M^{2n}$.

Let a mapping $q=\bar{q}(\mu ;h)$,
$(\mu ;h)\in M_{h,\tau }^{2n}$, be an
appropriate algebraic solution to equations (2.34). Substituting this into
the characteristic equations (2.32), one obtains, owing to (2.34), the
following characteristic Picard-Fuchs type [18,~19]  equations on
 $T^{\ast }(M_{h}^{n})$:
\begin{equation}
\partial w_{i}(\mu _{i};h)/\partial h_{j}={\bf \bar{P}}_{ji}(\mu
_{i},w_{i};h)
\end{equation}
for all  $i,j=\overline{1,n}$, where by def\/inition,
\begin{equation}
\ba{l}
\ds {\bf \bar{P}}_{ji}(\mu _{i},w_{i};h)=\left. {\bf P}_{ji}(\mu ;w;h)\right|
_{M_{h}^{n}\simeq \otimes _{j}^{n}{\mathbb S}_{j}^{1}}
\vspace{3mm}\\
\ds \ds \qquad =\sum_{i=1}^{n}\bar{h}_{js}(q(\mu ;h),
w\partial \mu /\partial q(\mu ;h))\partial q_{s}/\partial \mu _{i}.
\ea
\ee
As a result of the above computations one can formulate the following main
theorem.

\medskip

\noindent
{\bf Theorem 2.4.} {\it The imbedding (2.5) for the integral submanifold
$M_{h}^{n}\subset M^{2n}$ (compact and connected), parametrized
by a regular parameter $h\in {\mathcal G}^{\ast },$ is an algebraic
solution (up to diffeomorphism) to the set of characteristic Picard-Fuchs
type equations (2.35) on $T^{\ast }(M_{h}^{n})$, and can be
represented in general case~[19] in the following algebraic-geometric
form:
\begin{equation}
w_{j}^{n_{j}}+\sum_{s=1}^{n}c_{js}(\lambda ;h)w_{j}^{n_{j}-s}=0,
\end{equation}
where $c_{js}:{\mathbb R}\times {\mathcal G}^{\ast }\rightarrow
{\mathbb R}$, $s,j=\overline{1,n}$ are algebraic expressions, depending only
on the functional structure of the original abelian Lie algebra
${\mathcal G}$ of invariants on $M^{2n}$.
In particular, if the right-hand side of
the characteristic equations (3.5) is independent of
$h\in {\mathcal G}^{\ast }$, then this dependence will be linear in
$h\in {\mathcal G}^{\ast }$.}

\medskip

It should be noted here that some ten years ago an attempt was made in
[18,~19] to describe the explicit algebraic form of the Picard-Fuchs type
equations (2.35) by means of straightforward calculations for the well known
completely integrable Kowalewskaya top Hamiltonian system. The idea
suggested in [18,~19] was in some aspects very close to that devised
independently and thoroughly analyzed in~[11] which did not consider the
explicit form of the algebraic curves (2.37) starting from an abelian Lie
algebra  ${\mathcal G}$  of invariants  on a canonically symplectic phase
space $M^{2n}$.

As is well-known, a set of algebraic curves (2.31), prescribed via the above
algorithm, to a given apriori abelian Lie algebra ${\mathcal G}$ of
invariants on  the canonically symplectic phase space $M^{2n}=T^{\ast }
({\mathbb R}^{n})$ can be realized by means of a set of $n_{j}$-sheeted
Riemannian surfaces $\Gamma _{h}^{n_{j}}$, $j=\overline{1,n}$, covering
the corresponding real-valued cycles ${\mathbb S}_{j}^{1}$,
$j=\overline{1,n}$, which generate the corresponding homology
group $H_{1}({\mathbb T}^{n};{\mathbb Z})$ of the Arnold torus
${\mathbb T}^{n}\simeq \otimes _{j=1}^{n}{\mathbb S}_{j}^{1}$
dif\/feomorphic to the integral submanifold  $M_{h}^{n}\subset
M^{2n}$.

Thus, upon solving the set of algebraic equations (2.37) with respect to
functions $ w_{j}:{\mathbb S}_{j}^{1}\times {\mathcal G}^{\ast }\rightarrow
{\mathbb R}$,  $j=\overline{1,n}$, from (2.29) one obtains a vector parameter
$\tau =(t_{1},\ldots,t_{n})\in {\mathbb R}^{n}$ on
$M_{h}^{n}$ explicitly
described by means of the following Abelian type equations:
\begin{equation}
t_{j}=\sum_{s=1}^{n}\int_{\mu _{s}^{0}}^{\mu _{s}}d\lambda
\partial w_{s}(\lambda ;h)/\partial h_{j}
=\sum_{s=1}^{n}\int_{\mu _{s}^{0}}^{\mu _{s}}d\lambda
{\bf \bar{P}}_{js}(\lambda ,w_{s};h),
\ee
where $j=\overline{1,n}$, $(\mu ^{0};h)\in{ \otimes_{j=1}^n}
\Gamma _{h}^{n_{j}}\times
{\mathcal G}^{\ast }$. Using the expression (2.28) and recalling that the
generating function  $S:M_{h}^{n}\times {\mathbb R}^{n}\rightarrow
{\mathbb R}$ is a one-valued
mapping on an appropriate covering space
$\left(\bar{M}_{h}^{n};H_{1}(M_{h}^{n};{\mathbb Z})\right)$,
one can construct via the method of
Arnold~[4] the so called action-angle coordinates on $M_{h}^{n}$. Denote the
basic oriented cycles on $M_{h}^{n}$ by $\sigma _{j}\subset M_{h}^{n}$,
 $j=\overline{1,n}$.
These cycles together with their duals generate homology
group $H_{1}(M_{h}^{n};{\mathbb Z)}\simeq H_{1}({\mathbb T}^{n};
{\mathbb Z}) =\oplus _{j=1}^{n}{\mathbb Z}_{j}$.
In virtue of  the dif\/feomorphism $M_{h}^{n}\simeq \otimes _{j=1}^{n}
{\mathbb S}_{j}^{1}$ described above,  there
is a one-to-one correspondence between the basic cycles of
$H_{1}(M_{h}^{n};{\mathbb Z)}$ and those on the algebraic curves $\Gamma
_{h}^{n_{j}}$, $j=\overline{1,n}$, given by (2.37):
\begin{equation}
\rho :H_{1}(M_{h}^{n};{\mathbb Z})\rightarrow \oplus_{j=1}^{n}
{\mathbb Z}_{j}\sigma _{h,j},
\end{equation}
where $\sigma _{h,j}\subset \Gamma _{h}^{n_{j}}$, $j=\overline{1,n}$ are
the corresponding real-valued cycles on the Riemann surfaces $\Gamma
_{h}^{n_{j}}$, $j=\overline{1,n}$.

Assume that the following meanings of the mapping (2.39)  are prescribed:
\begin{equation}
\rho (\sigma _{i}):=\oplus _{j=1}^{n}n_{ij}\sigma _{h,j}
\end{equation}
for each  $i=\overline{1,n}$, where
$n_{ij}\in {\mathbb Z}$, $i,j=\overline{1,n}$
some f\/ixed integers. Then following the Arnold
construction [4,~18], one obtains the set of so called action-variables on
$M_{h}^{n}\subset M^{2n}$:
\begin{equation}
\gamma _{j:}:=\frac{1}{2\pi }\oint_{\sigma
_{j}}dS=\sum_{s=1}^{n}n_{js}\oint_{\sigma _{h,s}}d\lambda
w_{s}(\lambda ;h),
\end{equation}
where  $j=\overline{1,n}$. It is easy to show [4,~16], that expressions
(2.41) naturally def\/ine an a.e. dif\/ferentiable invertible mapping
\begin{equation}
\xi :{\cal G}^{\ast }\rightarrow {\mathbb R}^{n},
\end{equation}
which enables one to treat the integral submanifold $M_{h}^{n}$ as a
submanifold $M_{\gamma }^{n}$ $\subset M^{2n},$ where
\begin{equation}
M_{\gamma }^{n}:=\left\{u\in M^{2n}:\xi (h)=\gamma \in {\mathbb R}^{n}
\right\}.
\end{equation}
But, as was demonstrated in [18, 32], the functions (2.43) do not in general
generate a global foliation of the phase space $M^{2n}$, as they are
connected with both topological and analytical constraints. Since the
functions (2.41) are evidently also commuting invariants on $M^{2n}$, one
can def\/ine a further canonical transformation of the phase space
$M^{2n}$, generated by the following relationship on $M_{h,\tau }^{2n}$:
\begin{equation}
\sum_{j=1}^{n}w_{j} d\mu _{j}+\sum_{j=1}^{n}\varphi _{j}
d\gamma _{j}=dS(\mu ;\gamma ),
\end{equation}
where $\varphi =(\varphi _{1},\ldots,\varphi _{n})\in {\mathbb T}^{n}$
are the so called angle-variables on the torus
${\mathbb T}^{n}\simeq M_{h}^{n}$  and
 $S:M_{\gamma }^{n}\times {\mathbb R}^{n}\rightarrow {\mathbb R}$ is the
corresponding generating function. Whence it follows easily from (2.28) and
(2.38) that
\begin{equation}
\varphi _{j} := \frac{\partial S(\mu ;\gamma )}{\partial \gamma_{j}}
=\sum_{s=1}^{n}\frac{\partial S(\mu ;\gamma (h))}{\partial h_{s}}
\frac{\partial h_{s}}{\partial \gamma _{j}}
 =\sum_{s=1}^{n}t_{s}\omega _{sj}(\gamma ),
\qquad \frac{1}{2\pi }\oint_{\sigma _{j}}d\varphi _{k} =\delta _{jk},
\ee
where $\Omega :=\left\{\omega _{sj}:{\mathbb R}^{n}\rightarrow {\mathbb R},\
s,j=\overline{1,n}\, \right\}$ is the so called~[4] frequency matrix, which is
a.e. invertible on the integral submanifold $M_{\gamma }^{n}\subset
M^{2n}$. As an evident result of (2.45), we claim that the evolution of any
vector f\/ield $K_{a}\in \Gamma (M^{2n})$ for $a\in {\mathcal G}$  on the
integral submanifold $M_{\gamma }^{n}\subset M^{2n}$ is quasiperiodic
with a set of frequencies generated by the matrix
$\Omega \stackrel{\mbox{\scriptsize a.e.}}{\in }
\mbox{Aut}({\mathbb R}^{n})$, def\/ined above.
As examples showing the ef\/fectiveness of
 the above method of construction of integral submanifold
imbeddings for abelian integrable  Hamiltonian systems, one can verify
the Liouville-Arnold integrability of all Henon-Heiles and Neumann type
systems described in detail in [21,~22]; however, we shall not dwell on this
here.

\section{Integral submanifold imbedding problem for a nonabelian Lie algebra
of invariants}

\setcounter{equation}{0}

{\bf 3.1.} We shall assume below that  there is given a
Hamiltonian  vector f\/ield $K\in \Gamma (M^{2n})$ on the canonically
symplectic phase space  $M^{2n}=T^{\ast }({\mathbb R}^{n})$,
$n\in {\mathbb Z}_{+}$, which is endowed with
a nonabelian Lie algebra ${\mathcal G}$ of
invariants, satisfying all the conditions of the Mishchenko-Fomenko
Theorem~1.5, that is
\begin{equation}
\dim {\mathcal G}+\mbox{rank}\; {\mathcal G} =\dim M^{2n}.
\end{equation}

Then, as was proved above, an integral submanifold $M_{h}^{r}\subset
M^{2n}$ at a regular element $h\in {\mathcal G}^{\ast }$ is
$\mbox{rank}\; {\mathcal G}=r$-dimensional
and dif\/feomorphic (when compact and connected)  to the
standard $r$-dimensional torus ${\mathbb T}^{r}\simeq \otimes _{j=1}^{r}
{\mathbb S}_{j}^{1}$. It is natural to ask the following question:  How does
one construct the corresponding integral submanifold imbedding
\begin{equation}
\pi _{h}:M_{h}^{r}\rightarrow M^{2n},
\end{equation}
which characterizes all possible orbits of the dynamical system
$K\in\Gamma (M^{2n})$?

Having gained some experience in constructing the imbedding (3.2) in the
case of the abelian Liouville-Arnold theorem on integrability by
quadratures, we proceed below to study the integral submanifold
$M_{h}^{r}\subset M^{2n}$ by means of Cartan's theory [3, 12, 16, 22] of the
integrable ideals in the Grassmann algebra $\Lambda (M^{2n})$.
Let ${\mathcal I(G}^{\ast })$ be an ideal in
$\Lambda (M^{2n})$, generated by independent dif\/ferentials
$dH_{j}\in \Lambda ^{1}(M^{2n})$, $j=\overline{1,k}$,
on an open neighborhood $U(M_{h}^{r})$,
where by def\/inition, $r=\dim {\mathcal G}$.
The ideal  ${\mathcal I(G}^{\ast })$ is obviously Cartan integrable
[23,~16] with the integral submanifold $M_{h}^{r}\subset M^{2n}$ (at a
regular element $h\in {\mathcal G}^{\ast })$, on which it vanishes, that is
$\pi _{h}^{\ast }{\mathcal I(G}^{\ast })=0$. The dimension $\dim
M_{h}^{r}=\dim M^{2n}- \dim{\mathcal G}=r=\mbox{rank}\;{\mathcal G}$
due to the condition (3.1) imposed on the Lie algebra
${\mathcal G}$. It is useful to note here
that owing to the inequality $r\leq k$ for the rank ${\mathcal G}$, one
readily obtains from (3.1) that $\dim{\mathcal G}=k$ $\geq n$.
Since each base element $H_{j}\in {\mathcal G}$,
$j=\overline{1,k}$, generates a
symplectically dual vector f\/ield $K_{j}\in \Gamma (M^{2n})$,
$j=\overline{1,k}$, one can try to study the corresponding dif\/ferential system
$K({\mathcal G)}$ which is also Cartan integrable on the entire open
neighborhood $U(M_{h}^{r})\subset M^{2n}$. Denote the
corresponding dimension of the integral submanifold by $\dim M_{h}^{k}=\dim
K({\mathcal G)=}k$.
Consider now an abelian dif\/ferential system $K({\mathcal G}_{h})
\subset K({\mathcal G})$, generated by the Cartan subalgebra
${\mathcal G}_{h}\subset {\mathcal G}$ and its integral submanifold
$\bar{M}_{h}^{r}\subset U(M_{h}^{r})$. Since the Lie
subgroup $G_{h}=\exp {\mathcal G}_{h}$ acts on the integral
submanifold $M_{h}^{r}$ invariantly (see Chapter~1) and $\dim
\bar{M}_{h}^{r} = \mbox{rank}\; {\mathcal G}=r$, it
follows that $\bar{M}_{h}^{r}=M_{h}^{r}$. On the other hand, the system
$K({\mathcal G}_{h})\subset K({\mathcal G})$ by def\/inition, meaning that
the integral submanifold $M_{h}^{r}$ is an invariant part of the
integral submanifold $M_{h}^{k}\subset U(M_{h}^{r})$ with respect to the
Lie group $G=\exp {\mathcal G}$ -- action on $M_{h}^{k}$.
In this case one has the following result.

\medskip

\noindent
{\bf Lemma 3.1.} {\it There exist just $(n-r)\in {\mathbb Z}_{+}$
vector fields $\tilde{F}_{j}\in K({\mathcal G})/K({\mathcal G}_{h})$,
$j=\overline{1,n-r}$, for which
\begin{equation}
\omega ^{(2)}(\tilde{F}_{i},\tilde{F}_{j})=0
\end{equation}
on $U(M_{h}^{r})$ for all $i,j=\overline{1,n-r}$.}

\medskip

$\blacktriangleleft$ {\bfseries \itshape Proof.} It is obvious that the matrix
$\omega (\tilde{K}):=\left\{\omega ^{(2)}(\tilde{K}_{i},\tilde{K}_{j}):\
i,j=\overline{1,k}\, \right\}$ has on $U(M_{h}^{r})$ the
$\mbox{rank}\; \omega (\tilde{K})=k-r$, since
$\dim_{{\mathbb R}}\ker(\pi _{h}^{\ast }\omega ^{(2)})=
\dim _{{\mathbb R}}(\pi_{h\ast }K({\mathcal G}_{h}))=r$ on
$M_{h}^{r}$  at the regular element
$h\in {\mathcal G}^{\ast }$. Let us now complexify the tangent space
$T(U(M_{h}^{r}))$ using its even dimensionality. Whence  one can easily
deduce that on $U(M_{h}^{r})$  there exist just $(n-r)\in {\mathbb Z}_{+}$
vectors (not vector f\/ields!)  $\tilde{K}_{j}^{{\mathbb C}}\in K^{{\mathbb C}
}({\mathcal G})/K^{{\mathbb C}}({\mathcal G}_{h})$,
$j=\overline{1,n-r}$, from the complexif\/ied~[24] factor space
$K^{{\mathbb C}}({\mathcal G})/K^{{\mathbb C}}({\mathcal G}_{h})$.
To show this,  let us reduce the skew-symmetric matrix
$\omega (\tilde{K})\in \mbox{Hom}({\mathbb R}^{k-r})$  to its selfadjoint
equivalent $\omega (\tilde{K}^{{\mathbb C}})\in \mbox{Hom}({\mathbb C}^{n-r})$,
having taken into account that $\dim_{{\mathbb R}}
{\mathbb R}^{k-r}=\dim _{{\mathbb R}}{\mathbb R}^{k+r-2r}=
\dim _{{\mathbb R}}{\mathbb R}^{2(n-r)}=\dim _{{\mathbb C}
}{\mathbb C}^{n-r}$. Let now $f_{j}^{{\mathbb C}}\in {\mathbb C}^{n-r}$,
 $j=\overline{1,n-r}$,
be eigenvectors of the nondegenerate selfadjoint matrix
$\omega (\tilde{K}^{{\mathbb C}})\in \mbox{Hom}({\mathbb C}^{n-r})$, that is
\begin{equation}
\omega \left(\tilde{K}^{{\mathbb C}}\right) f_{j}^{{\mathbb C}}=\tilde{\lambda}_{j}
f_{j}^{{\mathbb C}},
\end{equation}
where $\tilde{\lambda}_{j}\in {\mathbb R}$,
$j=\overline{1,n-r}$, and for
all $i,j=\overline{1,n-r}$,  $\langle f_{i}^{{\mathbb C}},f_{j}^{{\mathbb C}
}\rangle =\delta _{i,j}$.
The above obviously means that in the basis $\left\{f_{j}^{{\mathbb C}}
\in K^{{\mathbb C}}({\mathcal G})/K^{{\mathbb C}}({\mathcal G}_{h}):
j=\overline{1,n-r}\, \right\}$ the matrix $\omega (\tilde{K}^{{\mathbb C}})
\in \mbox{Hom}({\mathbb C}^{n-r})$ is strictly diagonal and representable as
\begin{equation}
\omega (\tilde{K}^{{\mathbb C}})=\sum_{j=1}^{n-r}\tilde{\lambda}_{j}
f_{j}^{{\mathbb C}}\otimes _{{\mathbb C}}f_{j}^{{\mathbb C}},
\end{equation}
where $\otimes _{\mathbb C}$ is the usual Kronecker tensor product
of vectors from ${\mathbb C}^{n-r}$. Owing to the construction of the
complexif\/ied matrix $\omega \left(\tilde{K}^{{\mathbb C}}\right)\in
\mbox{Hom}({\mathbb C}^{n-r})$,
 one sees that the space $K^{{\mathbb C}}({\mathcal G})/K^{{\mathbb C}}
({\mathcal G}_{h})\simeq {\mathbb C}^{n-r}$  carries a K\"{a}hler
structure~[24] with respect to which the following expressions
\begin{equation}
\omega \left(\tilde{K}\right)=\mbox{Im}\;\omega (\tilde{K}),
\qquad \langle \cdot ,\cdot\rangle_{{\mathbb R}}=\mbox{Re}\;
\langle\cdot ,\cdot\rangle
\end{equation}
hold. Making use now of the representation (3.5) and expressions (3.6),
one can f\/ind vector f\/ields $\tilde{F}_{j}\in K({\mathcal G})
/K({\mathcal G}_{h})$, $j=\overline{1,n-r}$, such that
\begin{equation}
\omega \left(\tilde{F}\right)=\mbox{Im}\; \omega \left(\tilde{F}^{{\mathbb C}}\right)=J,
\end{equation}
holds  on  $U(M_{h}^{r})$, where $J\in \mbox{Sp}({\mathbb C}^{n-r})$ is the
standard symplectic matrix, satisfying the complex structure~[24] identity
$J^{2}=-I$. In virtue of the normalization conditions
 $\langle f_{j}^{{\mathbb C}},f_{j}^{{\mathbb C}}\rangle=
\delta _{i,j}$, for all $i,j=\overline{1,n-r}$,
one easily infers from (3.7) that $\omega ^{(2)}\left(\tilde{F}_{i},
\tilde{F}_{j}\right)=0$  for all  $i,j=\overline{1,n-r}$, where
by def\/inition
\begin{equation}
\tilde{F}_{j}:=\mbox{Re}\;\tilde{F}_{j}^{{\mathbb C}}
\end{equation}
for all  $j=\overline{1,n-r}$, and this proves the lemma.~$\blacktriangleright$

Assume now that  the Lie algebra ${\mathcal G}$ of invariants on $M^{2n}$
has been split into a direct sum of subspaces as
\begin{equation}
{\mathcal G}={\mathcal G}_{h}\oplus \widetilde{{\mathcal G}}_{h},
\end{equation}
where ${\mathcal G}_{h}$ is the Cartan subalgebra at a regular element
$h\in {\mathcal G}^{\ast }$ (being abelian) and
$\widetilde{{\mathcal G}}_{h}\simeq {\mathcal G}/{\mathcal G}_{h}$
is the corresponding  complement to  ${\mathcal G}_{h}$.
Denote a basis of  ${\mathcal G}_{h}$  as  $\left\{\bar{H}_{i}\in
{\mathcal G}_{h}:i=\overline{1,r}\, \right\}$, where $\dim {\mathcal G}_{h}=
\mbox{rank}\; {\mathcal G}=k\in {\mathbb Z}_{+}$,
and correspondingly, a basis of  $\widetilde{{\mathcal G}}_{h}$  as
$\left\{\tilde{H}_{j}\in \widetilde{{\mathcal G}}_{h}\simeq {\mathcal G}/
{\mathcal G}_{h}\!:\right.$ $\left.j=\overline{1,k-r}\, \right\}$. Then, owing to the results
of Chapter~1, the following relationships hold:
\begin{equation}
\{\bar{H}_{i},\bar{H}_{j}\}=0, \qquad
h(\{\bar{H}_{i},\tilde{H}_{s}\})=0
\end{equation}
on the open neighborhood $U(M_{h}^{r})\subset M^{2n}$
for all $i,j=\overline{1,r}$  and
$s=\overline{1,k-r}$. We have as yet had nothing to
say of expressions $h(\{\tilde{H}_{s},\tilde{H}_{m}\})$ for
$s,m=\overline{1,k-r}$. Making use of the representation~(3.8) for our vector
f\/ields (if they exist) $\tilde{F}_{j}\in K({\mathcal G})/K({\mathcal G}
_{h})$, $j=\overline{1,n-r}$, one can write down the following
expansion:
\begin{equation}
\tilde{F}_{i}=\sum_{j=1}^{k-r}c_{ji}(h)\tilde{K}_{j},
\end{equation}
where $i_{\tilde{K}_{j}}\omega ^{(2)}:=-d\tilde{H}_{j}$,
$c_{ji}:{\mathcal G}^{\ast }\rightarrow {\mathbb R}$,
$i=\overline{1,n-r}$,  $j=\overline{1,k-r}$,
are real-valued functions on ${\mathcal G}^{\ast }$, being def\/ined uniquely
as a result of  (3.11). Whence it clearly follows that there exist
invariants $\tilde{f}_{s}:U(M_{h}^{r})\rightarrow {\mathbb R}$,
$s=\overline{1,n-r}$,  such that
\begin{equation}
-i_{\tilde{F}_{s}}\omega ^{(2)}=\sum_{j=1}^{k-r}c_{js}(h)
d\tilde{H}_{j}:=d\tilde{f}_{s},
\end{equation}
where $\tilde{f}_{s}=\sum\limits_{j=1}^{k-r}c_{js}(h)\tilde{H}_{j}$,
 $ s=\overline{1,n-r}$, holds on $U(M_{h}^{r})$.

{\bf 3.2.} To proceed further, let us look at the following identity which is
similar to~(2.2):
\begin{equation}
\left(\otimes _{j=1}^{r}i_{\bar{K}_{j}}\right)\left(\otimes _{s=1}^{n-r}i_{_{\tilde{F}
_{s}}}\right)\left(\omega ^{(2)}\right)^{n+1}=0=\pm (n+1)!\left(\wedge _{j=1}^{r}d\bar{H}
_{j}\right)\left(\wedge _{s=1}^{n-r}d\tilde{f}_{s}\right)\wedge \omega ^{(2)},
\end{equation}
on  $U(M_{h}^{r})$. Whence, the following result is easily obtained using
Cartan theory~[3,~16]:

\medskip

\noindent
{\bf Lemma 3.2.} {\it The symplectic structure $\omega ^{(2)}\in \Lambda
^{2}(U(M_{h}^{r}))$  has the following canonical representation:
\begin{equation}
\left. \omega ^{(2)}\right| _{U(M_{h}^{r})}=\sum_{j=1}^{r}d\bar{H}_{j}\wedge
\bar{h}_{j}^{(1)}+\sum_{s=1}^{n-r}d\tilde{f}_{s}\wedge
\tilde{h}_{s}^{(1)},
\end{equation}
where $\bar{h}_{j}^{(1)}$, $\tilde{h}_{s}^{(1)}$
$\in \Lambda ^{1}(U(M_{h}^{r}))$, $j=\overline{1,r}$, $s=\overline{1,n-r}$.}

\medskip

The expression (3.14) obviously means, that on $U(M_{h}^{r})\subset M^{2n}$
the dif\/ferential 1-forms $\bar{h}_{j}^{(1)}$, $\tilde{h}_{s}^{(1)}
\in \Lambda ^{1}(U(M_{h}^{r}))$, $j=\overline{1,r}$, $s=\overline{1,n-r}$,
are independent together with exact  1-forms $d\bar{H}_{j}$, $j=\overline{1,r}$,
and $ d\tilde{f}_{s}$, $s=\overline{1,n-r}$. Since
$d\omega ^{(2)}=0$  on $M^{2n}$ identically, from (3.14)
one obtains that the dif\/ferentials  $d\bar{h}_{j}^{(1)}$,
$d\tilde{h}_{s}^{(1)} \in \Lambda ^{2}(U(M_{h}^{r}))$, $j=\overline{1,r}$,
$s=\overline{1,n-r}$, belong to the ideal
${\mathcal I}( \widetilde{{\mathcal G}}_{h})\subset {\mathcal I}({\mathcal G}^{\ast })$,
generated by exact forms $d\tilde{f}_{s}$, $s=\overline{1,n-r}$,
and $d\bar{H}_{j}$, $j=\overline{1,r}$, for all
regular $h\in {\mathcal G}^{\ast }$. Consequently, one obtains the following
analog of the Galisau-Reeb Theorem~2.1.

\medskip

\noindent
{\bf Theorem 3.3.} {\it Let a Lie algebra ${\mathcal G}$  of invariants
on the symplectic space $M^{2n}$ be nonabelian and satisfy the
Mishchenko-Fomenko condition~(3.1). At a regular element $h\in {\mathcal G}^{\ast }$
on some open  neighborhood $U(M_{h}^{r})$ of the
integral submanifold $M_{h}^{r}\subset M^{2n}$ there exist
differential 1-forms $\bar{h}_{j}^{(1)}$, $j=\overline{1,n}$,
and $\tilde{h}_{s}^{(1)}$, $s=\overline{1,n-r}$, satisfying the following properties:

i) $\left. \omega ^{(2)}\right|
_{U(M_{h}^{r})}=\sum\limits_{j=1}^{r}d\bar{H}_{j}\wedge \bar{h}
_{j}^{(1)}+\sum\limits_{s=1}^{n-r}d\tilde{f}_{s}\wedge \tilde{h}
_{s}^{(1)}$,
where $\bar{H}_{j}\in {\mathcal G}$, $j=\overline{1,r}$,
is a basis of the Cartan subalgebra ${\mathcal G}_{h}\subset {\mathcal G}$
 (being abelian), and $\tilde{f}_{s} \in {\mathcal G}$,  $s=\overline{1,n-r}$,
are  invariants from the complementary space $\widetilde{{\mathcal G}}_{h}
\simeq {\mathcal G}/{\mathcal G}_{h};$

ii) 1-forms $\bar{h}_{j}^{(1)} \in \Lambda^{1}(U(M_{h}^{r}))$,
$j=\overline{1,r}$, and $\tilde{h}_{s}^{(1)} \in \Lambda ^{1}(U(M_{h}^{r}))$,
$s=\overline{1,n-r}$, are exact on $M_{h}^{r}$ and
satisfy the equations: $\bar{h}_{j}^{(1)}(\bar{K}_{i})=\delta _{i,j}$ for all
$i,j=\overline{1,r}$, $\bar{h}_{j}^{(1)}(\tilde{F}_{s})=0$
and $\tilde{h}_{s}^{(1)}(\bar{K}_{j})=0$ for all $j=\overline{1,r}$,
$s=\overline{1,n-r}$, and $\tilde{h}_{s}^{(1)}(\tilde{F}_{m})=\delta _{s,m}$
for all $s,m=\overline{1,n-r}$.}

\medskip

$\blacktriangleleft$ {\bfseries \itshape Proof.} Obviously we need to prove only
the last statement {\it ii)}. Making use of Theorem~3.3, one f\/inds on the integral
submanifold $M_{h}^{r}\subset M^{2n}$ the dif\/ferential 2-forms
$d\bar{h}_{j}^{(1)}\in \Lambda ^{2}(U(M_{h}^{r}))$, $j=\overline{1,r}$,
and $d\tilde{h}_{s}^{(1)}\in \Lambda ^{2}(U(M_{h}^{r}))$, $s=\overline{1,n-r}$,
are identically  va\-ni\-shing. This means in
particular, owing  to the classical Poincar\'{e} lemma [1, 4, 16], that there
exist some exact 1-forms $d\bar{t}_{h,j}\in \Lambda ^{1}(U(M_{h}^{r}))$,
$j=\overline{1,r}$, and $d\tilde{t}_{h,s}\in \Lambda^{1}(U(M_{h}^{r}))$,
$s=\overline{1,n-r}$, where $\bar{t}_{h,j}: M_{h}^{r}\rightarrow {\mathbb R}$,
$j=\overline{1,r}$, and $\tilde{t}_{h,s}:M_{h}^{r}\rightarrow {\mathbb R}$,
$s=\overline{1,n-r}$, are smooth independent  a.e. functions on
$M_{h}^{r}$; they are one-valued on an appropriate covering of the manifold
$M_{h}^{r}\subset M^{2n}$ and supply global coordinates on the
integral submanifold $M_{h}^{r}$. Using the representation (3.14), one can
easily obtain that
\begin{equation}
-\left. i_{\bar{K}_{i}}\omega ^{(2)}\right|
_{U(M_{h}^{r})}=\sum_{j=1}^{r}d\bar{H}_{j}\bar{h}_{j}^{(1)}(\bar{K}_{i})
+\sum_{s=1}^{n-r}d\tilde{f}_{s}\tilde{h}_{s}^{(1)}(\bar{K}_{i})=d\bar{H}_{i}
\end{equation}
for all $i=\overline{1,r}$ and
\begin{equation}
-\left. i_{\tilde{F}_{m}}\omega ^{(2)}\right| _{U(M_{h}^{r})}=\sum_{j=1}^{r}d
\bar{H}_{j}\bar{h}_{j}^{(1)}(\tilde{F}_{m})+\sum_{s=1}^{n-r}d\tilde{f}_{s}
\tilde{h}_{s}^{(1)}(\tilde{F}_{m})=d\tilde{f}_{m}
\end{equation}
for all $m=\overline{1,n-r}$. Whence, from (3.15) it follows on that
on $U(M_{h}^{r})$,
\begin{equation}
\bar{h}_{j}^{(1)}(\bar{K}_{i})=\delta _{i,j}, \qquad \tilde{h}_{s}^{(1)}(\bar{K}_{i})=0
\end{equation}
for all $i,j=\overline{1,r}$ and  $s=\overline{1,n-r}$, and similarly,
from (3.16) it follows that on $U(M_{h}^{r})$,
\begin{equation}
\bar{h}_{j}^{(1)}(\tilde{F}_{m})=0, \qquad \tilde{h}_{s}^{(1)}(\tilde{F}_{m})=0
\end{equation}
for all  $j=\overline{1,r}$  and  $s,m=\overline{1,n-r}$. Thus the
theorem is proved.~$\blacktriangleright$

Having now def\/ined global evolution parameters  $t_{j}:M^{2n}\rightarrow
{\mathbb R}$,  $j=\overline{1,r}$, of the corresponding vector f\/ields
$\bar{K}_{j}=d/dt_{j}$, $j=\overline{1,r}$,  and local evolution
parameters  $\tilde{t}_{s}:M^{2n}\cap U(M_{h}^{r})\rightarrow {\mathbb R}$,
$s=\overline{1,n-r}$,  of the corresponding vector f\/ields $\left. \tilde{F}_{s}
\right| _{U(M_{h}^{r})}:=d/d\tilde{t}_{s}$,  $s=\overline{1,n-r}$,
one can easily see from (3.18) that the equalities
\begin{equation}
\left. t_{j}\right| _{U(M_{h}^{r})}=\bar{t}_{j}, \qquad
\left. \tilde{t}_{s}\right| _{U(M_{h}^{r})}=\tilde{t}_{h,s}
\end{equation}
hold for all  $j=\overline{1,r}$, $s=\overline{1,n-r}$,  up to
constant normalizations. Thereby, one can develop a new method, similar to
that of Chapter~2, for studying  the integral submanifold imbedding problem
in the case of the nonabelian Liouville-Arnold integrability theorem.

Before starting, it is interesting to note that the system of invariants
\[
\ {\mathcal G}_{\tau }:={\mathcal G}_{h}\oplus
\mbox{span}_{{\mathbb R}}\{\tilde{f}_{s}\in
{\mathcal G}/{\mathcal G}_{h}:s=\overline{1,n-r}\}
\]
constructed above, comprise a new  involutive (abelian) complete algebra
${\mathcal G}_{\tau }$, to which one can apply the abelian Liouville-Arnold
theorem on integrability by quadratures and the integral submanifold
imbedding theory devised in Chapter~2, in order to obtain exact solutions
by means of algebraic-analytical expressions. Namely, the following
corollary holds.

\medskip

\noindent
{\bf Corollary 3.5.} {\it Assume that a nonabelian Lie algebra ${\mathcal G}$
satisfies the Mishchenko-Fomenko condition (3.1) and $M_{h}^{r}\subset M^{2n}$
is its integral submanifold (compact and
connected) at a regular element $h\in {\mathcal G}^{\ast }$,  is
diffeomorphic to the standard torus ${\mathbb T}^{r}\simeq M_{h,\tau }^{n}$.
Assume also that  the dual complete  abelian algebra ${\mathcal G}_{\tau}$
$(\dim {\mathcal G}_{\tau }=n=1/2\dim M^{2n})$  of independent
invariants constructed above is globally defined. Then its integral
submanifold $M_{h,\tau }^{n}\subset M^{2n}$ is
diffeomorphic to  the standard torus  ${\mathbb T}^{n}\simeq M_{h,\tau }^{n}$,
and contains the torus ${\mathbb T}^{r}\simeq M_{h}^{r}$  as a
direct product with some completely degenerate torus ${\mathbb T}^{n-r}$,
that is $M_{h,\tau }^{n}\simeq M_{h}^{r}\times {\mathbb T}^{n-r}$.}

\medskip

Thus, having successfully applied the algorithm of Chapter 2 to the
algebraic-analytical characterization of integral submanifolds of a
nonabelian Liouville-Arnold integrable Lie algebra ${\mathcal G}$ of invariants
on the canonically symplectic manifold $M^{2n}\simeq T^{\ast }({\mathbb R}^{n})$,
one can produce a wide class of exact solutions represented by
quadratures -- which is just what we set out to f\/ind. At this point it is
necessary to note that up to now the (dual to ${\mathcal G}$) abelian complete
algebra ${\mathcal G}_{\tau }$ of invariants at a regular
$h\in {\mathcal G}^{\ast }$
was constructed only on some open neighborhood $U(M_{h}^{r})$
of the integral submanifold $M_{h}^{r}\subset M^{2n}$. As was mentioned
before, the global existence of the algebra ${\mathcal G}_{\tau }$ strongly
depends on the possibility of extending these invariants to the entire
manifold $M^{2n}$. This possibility is in 1--1 correspondence with the
existence of a global complex structure~[24] on the reduced integral
submanifold $\tilde{M}_{h,\tau }^{2(n-r)}:=M_{h}^{k}/G_{h}$, induced by the
reduced symplectic structure $\pi _{\tau }^{\ast }\omega ^{(2)}\in
\Lambda ^{2}(M_{h}^{k}/G_{h})$, where $\pi _{\tau }:M_{h}^{k}\rightarrow
M^{2n}$ is the imbedding for the integrable dif\/ferential system
$K({\mathcal G})\subset \Gamma (M^{2n})$, introduced above. If this is the case, the
resulting complexif\/ied manifold  $ ^{{\mathbb C}}\tilde{M}_{h,\tau
}^{n-r}\simeq \tilde{M}_{h,\tau }^{2(n-r)}$  will be endowed with a
K\"{a}hlerian structure, which makes it possible to produce the dual
abelian algebra ${\mathcal G}_{\tau }$ as a globally def\/ined set of invariants
on $M^{2n}$. This problem will be analyzed in more detail in Chapter~5.

\section{Examples}

\setcounter{equation}{0}

{\bf 4.1.} Below we consider some examples of  nonabelian
Liouville-Arnold integrability  by quadratures covered by Theorem~1.5.

\medskip

\noindent
{\bf Example 4.1.} {\bfseries \itshape Point vortices in the plane.}
Consider $n\in {\mathbb Z}_{+}$ point vortices on the plane  ${\mathbb R}^{2}$,
described by the Hamiltonian function
\begin{equation}
H=-\frac{1}{2\pi }\sum_{i\neq j=1}^{n}\xi _{i}\xi _{j}\ln \left\|
q_{i}-p_{j}\right\|
\end{equation}
with respect to the following partially canonical symplectic structure on
$M^{2n}\simeq T^{\ast }({\mathbb R}^{n})$:
\begin{equation}
\omega ^{(2)}=\sum_{j=1}^{n}\xi _{j}dp_{j}\wedge dq_{j},
\end{equation}
where $(p_{j},q_{j})\in {\mathbb R}^{2}$, $j=\overline{1,n}$, are
coordinates of the vortices in ${\mathbb R}^{2}$. There exist three additional
invariants
\begin{equation}
\ds P_{1}=\sum_{j=1}^{n}\xi _{j}q_{j}, \qquad P_{2}=\sum_{j=1}^{n}\xi
_{j}p_{j},
\qquad \ds P=\frac{1}{2}\sum_{j=1}^{n}\xi _{j}\left(q_{j}^{2}+p_{j}^{2}\right),
\ee
satisfying the following Poisson bracket conditions:
\begin{equation}
\ba{l}
\ds \{P_{1},P_{2}\}=-\sum_{j=1}^{n}\xi _{j}, \qquad \{P_{1},P\}=-P_{2},
\vspace{3mm}\\
\ds  \{P_{2},P\}=P_{1}, \qquad  \{P,H\}= 0 = \{P_{j},H\}.
\ea
\ee

It is evident, that invariants (4.1) and (4.3) comprise on
$\sum\limits_{j=1}^{n}\xi _{j}=0$ a four-dimensional Lie algebra
${\mathcal G}$, whose  $\mbox{rank}\; {\mathcal G}=2$.
Indeed, assume a regular vector $h\in {\mathcal G}^{\ast }$ is chosen,
and parametrized \ by real values $h_{j}\in {\mathbb R}$,
$j=\overline{1,4}$, where
\begin{equation}
h(P_{i})=h_{i}, \qquad  h(P)=h_{3}, \qquad  h(H)=h_{4},
\end{equation}
and $i=\overline{1,2}$. Then, one can easily verify that the element
\begin{equation}
Q_{h}=\left(\sum_{j=1}^{n}\xi _{j}\right)P-\sum_{i=1}^{n}h_{i}P_{i}
\end{equation}
belongs to the Cartan Lie subalgebra ${\mathcal G}_{h}\subset {\mathcal G}$,
that is
\begin{equation}
h(\{Q_{h},P_{i}\})=0, \qquad h(\{Q_{h},P\})=0.
\end{equation}

Since $\{Q_{h},H\}=0$ for all values $h\in {\mathcal G}^{\ast }$, we
claim that ${\mathcal G}_{h}=\mbox{span}_{{\mathbb R}}\{H,Q_{h}\}$ -- the Cartan
subalgebra of ${\mathcal G}$. Thus, $\mbox{rank}\; {\mathcal G}=
\dim {\mathcal G}_{h}=2$, and one comes right away that the condition  (3.1)
\begin{equation}
\dim M^{2n}=2n=\mbox{rank}\;{\mathcal G}+\dim {\mathcal G}=6
\end{equation}
holds only if $n=3$. Thereby, the following theorem is proved.

\medskip

\noindent
{\bf Theorem 4.1.} {\it The three -- vortex problem (4.1) on the plane
${\mathbb R}^{2}$  is nonabelian Liouville-Arnold integrable by quadratures
on the phase space $M^{6}\simeq T^{\ast }({\mathbb R}^{3})$  with the
symplectic structure (4.2).}

\medskip

As a result, the corresponding integral submanifold $M_{h}^{2}\subset M^{6}$
is two-dimensional and dif\/feomorphic (when compact
and connected) to the torus ${\mathbb T}^{2}\simeq M_{h}^{2}$, on
which the motions are quasiperiodic functions of the evolution parameter.

Concerning the Corollary 3.5, the dynamical system (4.1) is also abelian
Liouville-Arnold integrable with an extended integral submanifold
$M_{h,\tau }^{3}\subset M^{6}$, which can be found via the scheme suggested
in Chapter~3. Using simple calculations, one obtains an additional invariant
$Q=\left(\sum\limits_{j=1}^{3}\xi _{j}\right)P-\sum\limits_{i=1}^{3}P_{i}^{2}\notin
{\mathcal G}$, which commutes with $H$ and
$P$ of ${\mathcal G}_{h}$. Therefore, there
exists a new complete dual abelian algebra ${\mathcal G}_{\tau } =
\mbox{span}_{{\mathbb R}}
\{Q,P,H\}$  of  independent invariants  on  $M^{6}$ with
$\dim {\mathcal G}_{\tau }=3=1/2\dim M^{6}$, whose integral submanifold
$M_{h,\tau }^{3}\subset M^{6}$ (when compact and connected) is dif\/feomorphic
to the torus  ${\mathbb T}^{3}\simeq M_{h}^{2}\times {\mathbb S}^{1}$.

Note also here, that the above additional invariant  $Q\in {\mathcal G}_{\tau }$
can be naturally extended to the case of an arbitrary number  $n\in {\mathbb Z}_{+}$
of vortices as follows: $Q=\left(\sum\limits_{j=1}^{n}\xi_{j}\right)P-
\sum\limits_{i=1}^{n}P_{i}^{2}\in {\mathcal G}_{\tau }$, which obviously also
commutes with invariants (4.1) and (4.3) on the entire phase space $M^{2n}$.

\medskip

\noindent
{\bf Example 4.2.} {\bfseries \itshape A material point motion in a central  field.}
Consider the motion of a material point in the space ${\mathbb R}^{3}$
under a central potential f\/ield whose Hamiltonian
\begin{equation}
H=\frac{1}{2}\sum_{j=1}^{3}p_{j}^{2}+Q(\left\| q\right\| ),
\end{equation}
contains a central f\/ield $Q:{\mathbb R}_{+}\rightarrow {\mathbb R}$. The
motion is takes place in the canonical phase space $M^{6}=T^{\ast }
({\mathbb R}^{3})$, and possesses three additional invariants:
\begin{equation}
P_{1}=p_{2}q_{3}-pq,\qquad P_{2}=p_{3}q_{1}-p_{1}q_{3},\qquad
P_{3}=p_{1}q_{2}-p_{2}q_{1},
\end{equation}
satisfying the following Poisson bracket relations:
\begin{equation}
\{P_{1},P_{2}\}=P_{3}, \qquad \{P_{3},P_{1}\}=P_{2}, \qquad
\{P_{2},P_{3}\}=P_{1}.
\end{equation}
Since $\{H,P_{j}\}=0$  for all $j=\overline{1,3}$, one sees that
the problem under consideration has a four-dimensional Lie algebra
${\mathcal G}$ of invariants, isomorphic to the classical rotation Lie algebra
$so(3)\times {\mathbb R}\simeq {\mathcal G}$. Let us show that at a regular element
$h\in {\mathcal G}^{\ast }$  the Cartan subalgebra
${\mathcal G}_{h} \subset {\mathcal G}$  has the dimension $\dim {\mathcal G}_{h}=
2=\mbox{rank}\; {\mathcal G}$. Indeed, one easily verif\/ies that the invariant
\begin{equation}
P_{h}=\sum_{j=1}^{3}h_{j}P_{j}
\end{equation}
belongs to the Cartan subalgebra ${\mathcal G}_{h}$, that is
\begin{equation}
\{H,P_{h}\}=0, \qquad h(\{P_{h},P_{j}\})=0
\end{equation}
for all $j=\overline{1,3}$. Thus, as the Cartan subalgebra ${\mathcal G}_{h}
=\mbox{span}_{{\mathbb R}}\{H$ and $P_{h}\subset {\mathcal G}\}$, one gets
$\dim {\mathcal G}_{h} = 2=\mbox{rank}\; {\mathcal G}_{h}$,
and the Mishchenko-Fomenko  condition~3.1
\begin{equation}
\dim M^{6}=6=\mbox{rank}\; {\mathcal G}+\dim {\mathcal G}=4+2
\end{equation}
holds. Hence one can prove its integrability by quadratures  via the
nonabelian Liouville Liouville-Arnold Theorem~1.5 and obtain the following
theorem:

\medskip

\noindent
{\bf Theorem 4.3.} {\it It follows from Theorem 1.5 that the free material
point motion in ${\mathbb R}^{3}$ is a completely integrable by
quadratures dynamical system on the canonical symplectic phase space
$M^{6}=T^{\ast }({\mathbb R}^{3})$. The corresponding integral submanifold
$M_{h}^{2}\subset M^{6}$  at a regular element  $h\in {\mathcal G}^{\ast }$
(if compact and connected) is two-dimensional and
diffeomorphic to the standard torus ${\mathbb T}^{2}\simeq M_{h}^{2}$.}

\medskip

Making use of the integration algorithm devised in Chapters~1 and~2, one can
readily obtain the corresponding integral submanifold imbedding mapping
$\pi _{h}: M_{h}^{2}\rightarrow M^{6}$ by means of
algebraic-analytical expressions and transformations.

There are clearly many other interesting nonabelian Liouville-Arnold
integrable Hamiltonian systems on canonically symplectic phase spaces that
arise in applications, which can similarly be integrated using
algebraic-analytical means. We hope to study several of these systems in
detail elsewhere.

\section{Existence problem for a global set of invariants}

\setcounter{equation}{0}

{\bf 5.1} It was proved in Chapter 3, that
locally, in some open neighborhood $U(M_{h}^{r})\subset M^{2n}$ of the
integral submanifold  $M_{h}^{r}\subset M^{2n}$  one can f\/ind by
algebraic-analytical  means  just $n-r\in {\mathbb Z}_{+}$  independent
vector f\/ields $\tilde{F}_{j}\in K({\mathcal G})/K({\mathcal G}_{h})\cap
\Gamma (U(M_{h}^{r})$, $j=\overline{1,n-r}$, satisfying the condition~(3.3).
Since each vector f\/ield $\tilde{F}_{j}\in K({\mathcal G})/K({\mathcal G}_{h})$,
$j=\overline{1,n-r}$, is generated by an invariant  $\tilde{H}_{j}
 \in {\mathcal D}(U(M_{h}^{r}))$, $j=\overline{1,n-r}$, it follows
readily from (3.3) that
\begin{equation}
\{\tilde{H}_{i},\tilde{H}_{j}\}=0
\end{equation}
for all $i,j=\overline{1,n-r}$. Thus, on an open neighborhood
$U(M_{h}^{r})$ there exist just  $n-r$ invariants in addition to
$\tilde{H}_{j} \in {\mathcal D}(U(M_{h}^{r}))$, $j=\overline{1,n-r}$,
all of which are in involution. Denote as before this new set of invariants
as ${\mathcal G}_{\tau }$, keeping in mind that
$\dim {\mathcal G}_{\tau }=r+(n-r)=n\in {\mathbb Z}_{+}$.
Whence, on an open neighborhood  $U(M_{h}^{r})\subset M^{2n}$
we have constructed the set  ${\mathcal G}_{\tau }$ of just
$n=1/2\dim M^{2n}$ invariants commuting with each
other, thereby guaranteeing via the abelian Liouville-Arnold theorem its
local complete integrability by quadratures. Consequently, there exists
locally a mapping $\pi _{\tau }:M_{h,\tau }^{k}\rightarrow M^{2n}$,
where $M_{h,\tau }^{k}:=U(M_{h}^{r})\cap M_{\tau }^{k}$ is the integral
submanifold of the dif\/ferential system $K({\mathcal G})$,  and one can
therefore describe the behavior of integrable vector f\/ields on the reduced
manifold $\bar{M}_{h,\tau }^{2(n-r)} := M_{h,\tau }^{k-r}/G_{h}$.
\ For global integrability properties of a given set  ${\mathcal G}$ of
invariants on $\left(M^{2n},\omega ^{(2)}\right)$, satisfying the
Mishchenko-Fomenko condition~(3.1), it is necessary to have the additional
set of invariants $\tilde{H}_{j} \in {\mathcal D}(U(M_{h}^{r}))$,
$j=\overline{1,n-r}$, extended from $U(M_{h}^{r})$  to the entire
phase space $M^{2n}$. This problem evidently depends on  the existence of
extensions of vector f\/ields $\tilde{F}_{j}\in \Gamma (U(M_{h}^{r}))$,
$j=\overline{1,n-r}$, from the neighborhood
$U(M_{h}^{r})\subset M^{2n}$ to the whole phase space $M^{2n}$.
On the other hand, as stated
before, the existence of such a continuation depends intimately on the
properties of the complexif\/ied dif\/ferential system
$K^{{\mathbb C}}({\mathcal G})/K^{{\mathbb C}}({\mathcal G}_{h})$,
which has a nondegenerate complex  metric $\omega \left(\tilde{K}^{{\mathbb C}}\right):
T\left(\bar{M}_{h,\tau }^{2(n-r)}\right)^{{\mathbb C}}\times
T\left(\bar{M}_{h,\tau }^{2(n-r)}\right)^{{\mathbb C}}
\rightarrow {\mathbb C}$,
induced by the symplectic structure  $\omega ^{(2)}\in \Lambda^{2}(M^{2n}) $.
This point can be clarif\/ied more by using the notion [24--27] of a
K\"{a}hler manifold and some of the associated constructions
presented above. Namely, consider the local isomorphism  $T\left(\bar{M}
_{h,\tau }^{2(n-r)}\right)^{{\mathbb C}}\simeq T\left(^{{\mathbb C}}\bar{M}_{h,\tau }^{n-r}
\right)$,
where $^{{\mathbb C}}\bar{M}_{h,\tau }^{n-r}$ is the complex  $(n-r)$-dimensional
local integral submanifold of the complexif\/ied dif\/ferential
system $K^{{\mathbb C}}({\mathcal G})/K^{{\mathbb C}}({\mathcal G}_{h})$. This means
that the space $T\left(\bar{M}_{h,\tau }^{2(n-r)}\right)$ is endowed with the standard
almost complex structure
\begin{equation}
J:T\left(\bar{M}_{h,\tau }^{2(n-r)}\right)\rightarrow
T\left(\bar{M}_{h,\tau }^{2(n-r)}\right),
\qquad J^{2}=-1,
\end{equation}
such that the 2-form  $\omega (\tilde{K}):=\mbox{Im}\;
\omega \left(\tilde{K}^{{\mathbb C}}\right)
\in \Lambda ^{2}\left(\bar{M}_{h,\tau }^{2(n-r)}\right)$  induced from the
above metric on $T\left(^{{\mathbb C}}\bar{M}_{h,\tau }^{n-r}\right)$ is closed, that
is $d\omega (\tilde{K})=0$. If this is the case, the almost complex
structure on the manifold $T\left(\bar{M}_{h,\tau }^{2(n-r)}\right)$  is said to
be integrable. Def\/ine  the proper complex  manifold
$^{{\mathbb C}}\bar{M}_{h,\tau }^{n-r}$, on which one can then def\/ine global vector
f\/ields  $\tilde{F}_{j}\in K({\mathcal G})/K({\mathcal G}_{h})$, $j=\overline{1,n-r}$,
which are being sought for the involutive algebra ${\mathcal G}_{\tau }$ of
invariants on $M^{2n}$ to be integrable by quadratures via the abelian
Liouville-Arnold theorem. Thus the following theorem can be obtained.

\medskip

\noindent
{\bf Theorem 5.1.} {\it A nonabelian set ${\mathcal G}$  of
invariants on the symplectic space $M^{2n}\simeq T^{\ast }(R^{n})$,
satisfying the Mishchenko-Fomenko condition~3.1, admits
algebraic-analytical integration by quadratures for the integral
submanifold imbedding $\pi _{h}:M_{h}^{r}\rightarrow M^{2n}$,
if the corresponding complexified reduced manifold
$^{{\mathbb C}}\bar{M}_{h,\tau }^{n-r}\simeq \bar{M}_{h,\tau}^{2(n-r)}=
M_{h,\tau }^{k-r}/G_{h}$ of the differential system $K^{{\mathbb C}}({\mathcal G})
/K^{{\mathbb C}}({\mathcal G}_{h})$ is K\"{a}hlerian
with respect to the standard almost complex structure~(5.1) and the
nondegenerate complex  metric $ \omega \left(\tilde{K}^{{\mathbb C}}\right):
T\left(\bar{M}_{h,\tau }^{2(n-r)}\right)^{{\mathbb C}}\times
T\left(\bar{M}_{h,\tau }^{2(n-r)}\right)^{{\mathbb C}
}\rightarrow {\mathbb C}$,  induced by the symplectic structure
$\omega ^{(2)}\in \Lambda ^{2}(M^{2n})$ is integrable, that is $d\, \mbox{\rm Im}\;
\omega (\tilde{K}^{{\mathbb C}})=0$.}

\medskip

Theorem 5.1 shows, in particular, that nonabelian Liouville-Arnold
integrability by quadratures does not in general  imply integrability via
the abelian Liouville-Arnold theorem; it actually depends on certain
topological obstructions associated with the Lie algebra structure of
invariants ${\mathcal G}$  on the phase space $M^{2n}$. We hope to explore this
intriguing problem in another place.

\medskip


\subsection*{Acknowledgments}

The author is cordially thankful to Prof. Boris A. Kupershmidt (Space
Institute, University of Tennessee, Tullahoma USA) for many important
suggestions and critical comments concerning the exposition of the article.
He is also indebted to Profs. Denis Blackmore (Dept. of Mathem. Sciences at
NJIT, USA), Stefan Rauch-Wojciechowski (Dept. of Mathematics at Link\"oping
University, Sweden), Anatoliy M. Samoilenko (Institute of Mathematics at the
NAS, Kyiv, Ukraine), Andrzej Pelczar and Jerzy Ombach (Institute of
Mathematics at Jagiellonian University, Krakow, Poland), for many
stimulating discussions of the problems treated in this article. The author
is also thankful to the referee for very useful remarks and corrections. The
work on this article was in part supported by a local AGH-grant from
Krak\'{o}w.

\label{prykarp-lp}

\end{document}